\documentclass[12pt,a4paper,reqno]{amsart}
\usepackage[utf8]{inputenc}
\usepackage[T1]{fontenc}
\usepackage{amsmath}
\usepackage{amsthm}
\usepackage{amssymb}
\usepackage{stmaryrd}
\usepackage[abbrev]{amsrefs}
\usepackage[dvipsnames]{xcolor}
\usepackage{bm}
\usepackage{bbm}
\usepackage{enumitem}
\usepackage{hyperref}
\usepackage[all]{xy}
\AtBeginDocument{\def\MR#1{}}
\makeatletter
\@namedef{subjclassname@2020}{%
\textup{2020} Mathematics Subject Classification}
\makeatother
\numberwithin{equation}{section}
\allowdisplaybreaks
\newtheorem{Question}{}

\newtheorem{question}[Question]{Question}
\newtheorem{Claim}{}

\newtheorem{claim}[Claim]{Claim}
\newtheorem{theorem}{Theorem}[section]
\newtheorem{lemma}[theorem]{Lemma}
\newtheorem{proposition}[theorem]{Proposition}
\newtheorem{corollary}[theorem]{Corollary}
\theoremstyle{remark}
\newtheorem{remark}[theorem]{Remark}

\newcommand{\cc}{\ensuremath{c}}
\newcommand{\st}{\ensuremath{\bm{s}}}
\newcommand{\ZZ}{\ensuremath{\mathbb{Z}}}
\newcommand{\NN}{\ensuremath{\mathbb{N}}}
\newcommand{\Sym}{\ensuremath{\mathbb{S}}}
\newcommand{\AB}{\ensuremath{\mathbb{A}}}
\newcommand{\BB}{\ensuremath{\mathbb{B}}}
\newcommand{\RR}{\ensuremath{\mathbb{R}}}

\newcommand{\FF}{\ensuremath{\mathbb{F}}}
\newcommand{\XX}{\ensuremath{{\mathbb{X}}}}
\newcommand{\YY}{\ensuremath{\mathbb{Y}}}
\newcommand{\UU}{\ensuremath{\mathbb{U}}}
\newcommand{\EE}{\ensuremath{\mathbb{E}}}
\newcommand{\LL}{\ensuremath{\mathbb{L}}}
\newcommand{\EB}{\ensuremath{\mathcal{E}}}
\newcommand{\It}{\ensuremath{\mathcal{I}}}
\newcommand{\Jt}{\ensuremath{\mathcal{J}}}

\newcommand{\Ft}{\ensuremath{\mathcal{F}}}

\newcommand{\Nt}{\ensuremath{\mathcal{N}}}
\newcommand{\YB}{\ensuremath{\mathcal{Y}}}
\newcommand{\XB}{\ensuremath{\mathcal{X}}}
\newcommand{\xx}{\ensuremath{\bm{x}}}
\newcommand{\ee}{\ensuremath{\bm{e}}}
\newcommand{\Gr}{\ensuremath{\mathcal{G}}}

\newcommand{\Ind}{\ensuremath{\mathbbm{1}}}
\newcommand{\lsdf}{\ensuremath{\bm{\varphi_l}}}
\newcommand{\usdf}{\ensuremath{\bm{\varphi_u}}}
\DeclareMathOperator*{\spn}{span}

\DeclareMathOperator{\supp}{supp}
\DeclareMathOperator{\sgn}{sign}
\newcommand{\abs}[1]{\left\lvert#1\right\rvert}
\newcommand{\norm}[1]{\left\lVert#1\right\rVert}
\newcommand{\enfloor}[1]{\left\lfloor#1\right\rfloor}

\newcommand{\enbrace}[1]{\left\lbrace#1\right\rbrace}
\newcommand{\enbrak}[1]{\left[#1\right]}
\newcommand{\enpar}[1]{\left(#1\right)}
\newcommand{\uii}[1]{\left[#1\ldotp\ldotp\nobreak\infty\right)}
\newcommand{\bii}[2]{\left[#1\ldotp\ldotp\nobreak#2\right]}
\begin{document}
% ------------------------------------------------------------------------
\title[Counterexamples for the TGA in Quasi-Banach Spaces]{When Greedy Approximation Breaks: Counterexamples in Quasi-Banach Spaces}
% ------------------------------------------------------------------------
\author[F. Albiac]{Fernando Albiac}
\address{Department of Mathematics, Statistics, and Computer Sciences--Ina\-Mat$^2$ \\
Universidad P\'ublica de Navarra\\
Campus de Arrosad\'{i}a\\
Pamplona\\
31006 Spain}
\email{fernando.albiac@unavarra.es}
% ------------------------------------------------------------------------
\author[J. L. Ansorena]{Jos\'e L. Ansorena}
\address{Department of Mathematics and Computer Sciences\\
Universidad de La Rioja\\
Logro\~no\\
26004 Spain}
\email{joseluis.ansorena@unirioja.es}
% ------------------------------------------------------------------------
\author[M. Berasategui]{Miguel Berasategui}
\address{Miguel Berasategui\\ - Departamento de Matemática, - - Pab I, Facultad de Ciencias Exactas y Naturales\\
Universidad de Buenos Aires\\ Buenos Aires, 1428, Argentina}
\email{mberasategui@dm.uba.ar}
%--------------------------------------------------------------------------
\author[P. Bern\'a]{Pablo M. Bern\'a}
\address{Pablo M. Bern\'a\\
Departamento de Matem\'aticas, CUNEF Universidad\\
Madrid, 28040 Spain.}
\email{pablo.berna@cunef.edu}
%--------------------------------------------------------------------------
\begin{abstract}
We construct two counterexamples that resolve long-standing open problems on greedy approximation theory with respect to bases, posed in [F. Albiac et al., Dissertationes Math. 560 (2021)] and restated in [F.\@ Albiac, J.\@ L.\@ Ansorena, V.\@ Temlyakov, J.\@ Approx.\@ Theory 307 (2025)]. Our first result exhibits a quasi-Banach space $\XX$ with an almost greedy basis which, when transported to the Banach envelope of $\XX$, ceases to be quasi-greedy. This shows that the passage to the Banach envelope, although it preserves linear and lattice structure, may radically disrupt the performance of the thresholding greedy algorithm, to the extent that in some respects it could perform better in a quasi-Banach space than in its Banach envelope. Our second result constructs an almost greedy Markushevich basis in a nonlocally convex quasi-Banach space $\YY$ which fails to be a Schauder basis under any reordering. Together, these examples highlight that local convexity and the Banach envelope construction play an unexpectedly active role in shaping greedy approximation phenomena, revealing structural differences between Banach and quasi-Banach spaces that go beyond the classical theory of bases.
\end{abstract}
%--------------------------------------------------------------------------
\thanks{F.\@ Albiac and J.\@ L.\@ Ansorena acknowledge the support of the Spanish Ministry for Science and Innovation under Grant PID2022-138342NB-I00 for \emph{Functional Analysis Techniques in Approximation Theory and Applications (TAFPAA)}}
%--------------------------------------------------------------------------
\thanks{M.\@ Berasategui is supported by Grants CONICET PIP11220200101609CO y ANPCyT PICT 2018-04104 (Consejo Nacional de Investigaciones Científicas y Técnicas and Agencia Nacional de Promoción de la Investigación, el Desarrollo Tecnológico y la Innovación, Argentina)}
%--------------------------------------------------------------------------
\thanks{Pablo M.\@ Berná is supported by Grant PID2022-142202NB-I00 (Agencia Estatal de Investigación, Spain)}
%--------------------------------------------------------------------------
\maketitle
%-------------------------------------------------------- ------------------
\section{Introduction}\noindent
%--------------------------------------------------------------------------
The theory of greedy approximation investigates the performance of the thresholding greedy algorithm (TGA for short) with respect to bases in Banach and quasi-Banach spaces. At a technical level, the algorithm retains those basis coefficients of largest magnitude and discards the rest. At a conceptual level, it provides a lens through which one can investigate the interaction between the geometry of the ambient space and the structural properties of its bases.

Over the last two decades, this viewpoint has produced a rich body of results, including characterizations of new types of bases, such as quasi-greedy and almost greedy bases, and a systematic understanding of their connections with classical notions such as unconditionality and democracy.

The systematic study of greedy approximation with respect to bases in quasi-Banach spaces was initiated in the monograph \cite{AABW2021}, where Albiac, Ansorena, Bern\'a, and Wojtaszczyk extended the classical Banach space theory into the quasi-Banach framework. This extension is far from formal and requires having to reprove the foundational theorems with the invention of novel techniques because of the failure of the Hahn-Banach theorem and the duality theory when local convexity is lifted. Many spaces of central importance in analysis, such as $L_p$ for $0<p<1$ and Hardy or Besov spaces in the nonlocally convex regime, are naturally quasi-Banach, hence understanding the behaviour of greedy algorithms in these settings is therefore essential if the theory is to interact fully with harmonic analysis and approximation theory. Moreover, quasi-Banach spaces are intimately tied to their Banach envelopes, the closest Banach spaces into which they embed.

The relationship between a space and its envelope is often subtle and revealing, and the study of greedy approximation provides a new lens through which to explore this interaction. While in Banach spaces the theory is now mature, extending it to quasi-Banach spaces has revealed new phenomena and prompted a number of fundamental open problems. Two of these problems were posed in \cite{AABW2021} and restated in the recent survey \cite{AAT2025}.
\begin{question}[\cite{AABW2021}*{Problem 13.5}, cf.\@ \cite{AAT2025}*{Problem 16}]\label{qt:A}
Must a quasi-greedy basis in a quasi-Banach space remain quasi-greedy when considered in the Banach envelope?
\end{question}

\begin{question}[\cite{AABW2021}*{Problem 13.6}, cf.\@ \cite{AAT2025}*{Problem 14}]\label{qt:B}
([3, Problem 13.6], cf.\@ [2, Problem 14]) Can there exist a quasi-greedy basis that fails to be a Schauder basis under every possible reordering?
\end{question}

In this paper we give counterexamples that resolve both problems. Our first construction shows that a quasi-Banach space $\XX$ may carry an almost greedy Schauder basis whose image in the Banach envelope of $\XX$ is not quasi-greedy. This demonstrates that the Banach envelope, while designed to ``regularize'' a quasi-Banach space, can significantly distort the behaviour of bases with respect to greedy approximation. Although quasi-greedy bases always retain a flavour of unconditionality, the loss of linearity in the property takes its toll when passing to the Banach envelope. Our second construction provides an almost greedy Markushevich basis in a nonlocally convex quasi-Banach space $\YY$ which cannot be rearranged into a Schauder basis. This confirms that the success of the greedy algorithm is essentially independent of basis order, and that one should not expect Schauder-type structure to be present in general.

These results highlight two structural lessons for Banach space theory. First, local convexity is not merely a technical convenience: oddly enough, when it is removed, the TGA may actually perform better than in the corresponding Banach envelope! Second, the independence of greedy approximation from basis order makes the Schauder condition unnatural in this context, and our example shows definitively that greedy behaviour can survive in the absence of any Schauder ordering. In both directions, the counterexamples illustrate that the passage from Banach to quasi-Banach settings is not smooth, but instead uncovers fundamentally new interactions between geometry, convexity, and bases.

Throughout the paper we adhere to the standard notation of Banach space theory and approximation theory, as set out in the monograph \cite{AlbiacKalton2016} and the aforementioned article \cite{AABW2021}. To avoid ambiguity, we emphasize only a few conventions that may be less common in the literature. This will allow us to keep the exposition concise without departing from the established framework.

Given families $(s_j)_{j\in J}$ and $(t_j)_{j\in J}$ in $[0,\infty]$, the symbol $s_j\lesssim t_j$ for $j\in J$ means the existence of $C\in(0,\infty)$ such that $s_j\le C t_j$ for all $j\in J$. If $s_j\lesssim t_j$ and $t_j\lesssim s_j$ for $j\in J$, we say that $(s_j)_{j\in J}$ and $(t_j)_{j\in J}$ are equivalent, and we put $s_j\approx t_j$ for $j\in J$.

Given a set $\Nt$, $\enbrak{\Nt}$ will be the set consisting of all subsets of $\Nt$, and $\enbrak{\Nt}^{<\infty}$ the set consisting of all finite subsets of $\Nt$. Similarly, $\enbrak{\Nt}^{\le m}$ (resp., $\enbrak{\Nt}^{\ge m}$) denotes the set consisting of all finite subsets of $\Nt$ of cardinality at most (resp., at least) $m\in\NN$.

To denote integer intervals, given $k$, $m\in\ZZ$ with $k\le m$ we set
\[
\bii{k}{m}=\ZZ\cap[k,m], \quad \uii{k}=\ZZ\cap[k,\infty).
\]

We will use $\cc_{00}$ to denote the vector space consisting of all eventually null scalar-valued sequences, and $\cc_0$ to denote the space of all scalar-valued sequences converging to $0$. Other more specific notation will be introduced in context, as needed.
%--------------------------------------------------------------------------
\section{Terminology and background}\noindent
%--------------------------------------------------------------------------
A \emph{quasi-norm} on a vector space $\XX$ over the real or complex field $\FF$ is a map $\norm{ \cdot}\colon \XX\to [0, \infty)$ so that there is a constant $\kappa$ such that
\begin{enumerate}[label=(Q.\arabic*),leftmargin=*]
\item\label{it:Q1} $\norm{f}>0$ for all $f\not=0$,
\item\label{it:Q2} $\norm{t f}=\abs{t} \norm{f}$ for all $t\in \FF$ and all $f\in \XX$, and
\item\label{it:modofconc} $\norm{f+g}\le \kappa\enpar{ \norm{f}+\norm{g}}$ for all $f$, $g\in \XX$.
\end{enumerate}
A quasi-normed space is a vector space equipped with a quasi-norm. The smallest constant $\kappa$ in \ref{it:modofconc} will be called the \emph{modulus of concavity} of the norm and of the space. Set
\[
B_\XX=\enbrace{f\in\XX\colon\norm{ f} \le 1}, \quad
S_\XX=\enbrace{f\in\XX\colon\norm{ f} = 1}
\]
for the \emph{unit ball} and the \emph{unit sphere} of $\XX$, respectively. The dilations of $B_\XX$ form a basis of neighborhoods of the origin for a vector topology, and any quasi-normed space will considered equipped with this vector topology. A \emph{quasi-Banach space} is a complete quasi-normed space.

Let $0<p\le 1$. If $\norm{\cdot}\colon \XX\to[0,\infty)$ fulfils \ref{it:Q1}, \ref{it:Q2} and
\begin{enumerate}[label=(Q.\arabic*),leftmargin=*, resume]
\item $\norm{f+g}^p \le \norm{f}^p+\norm{g}^p$ for all $f$, $g\in \XX$.
\end{enumerate}
we say that $\norm{\cdot}$ is a \emph{$p$-norm} and that $\XX$ is a $p$-normed space. Any $p$-norm is a quasi-norm with modulus of concavity at most $2^{1/p-1}$. A $p$-Banach space is a complete $p$-normed space. A Banach space is a $p$-Banach space for $p=1$, or a quasi-Banach space whose modulus of concavity is one.

A subset $D$ of a vector space $\XX$ is \emph{absolutely $p$-convex} if
\[
\sum_{j\in J} t_j f_j\in D
\]
for all finite families $(f_j,t_j)_{j\in J}$ in $D\times\FF$ with $\sum_{j\in J} \abs{t_j}^p \le 1$. A quasi-normed space $\XX$ is said to be \emph{locally $p$-convex} if the origin has an absolutely $p$-convex neighborhood. Each of the two following conditions is equivalent to local $p$-convexity.

\begin{enumerate}[label=(L.\arabic*),leftmargin=*]
\item The quasi-norm $\norm{\cdot}$ of $\XX$ is equivalent to a $p$-norm.
\item\label{it:LpC:b} There is a constant $C$ such that
\[
\norm{\sum_{j\in J} f_j} \le C \enpar{\sum_{j\in J} \norm{f_j}^p}^{1/p}
\]
for any finite family $(f_j)_{j\in J}$ in $\XX$.
\end{enumerate}
Indeed, if \ref{it:LpC:b} holds, then there is a $p$-norm $\norm{\cdot}_0$ on $\XX$ with
\[
\norm{f}_0 \le \norm{f} \le C \norm{f}_0, \quad f\in\XX.
\]

By Aoki--Rolewicz Theorem \cites{Aoki1942,Rolewicz1957}, every quasi-Banach space is locally $p$-convex for some $p\in(0,1]$. Quantitatively, we have the following.

\begin{theorem}[see \cite{KPR1984}]\label{thm:AR}
Let $(\XX,\norm{\cdot})$ be a quasi-Banach space with modulus of concavity at most $\kappa$, and define $p\in(0,1]$ by $\kappa=2^{1/p-1}$. Then \ref{it:LpC:b} holds with $C=4^{1/p}$.
\end{theorem}

We point out that quasi-norms may be discontinuous maps relative to the topology they induce (see \cite{Hyers1939}). Since $p$-norms, $0<p\le 1$, are continuous, Theorem~\ref{thm:AR} gives a continuous quasi-norm $\norm{\cdot}_0$ with
\[
\norm{f}_0 \le \norm{f} \le 4 \kappa^2 \norm{f}_0, \quad f\in\XX,
\]
where $\kappa$ is the modulus of concavity of $\XX$. This equivalence allows estimates obtained for continuous quasi-norms to be transferred to arbitrary quasi-norms. Notwithstanding, applying the following lemma often yields sharper inequalities.

\begin{lemma}\label{lem:Approx}
Let $(\XX,\norm{\cdot})$ be a quasi-Banach space with modulus of concavity at most $\kappa$. Then,
\[
\frac{1}{\kappa} \limsup_{j\in\NN} \norm{f_j} \le \norm{f} \le \kappa \liminf_{j\in\NN} \norm{f_j}
\]
for every $f\in\XX$ and every sequence $(f_j)_{j=1}^\infty $ in $\XX$ converging to $f$.
\end{lemma}

\begin{proof}
The inequality
\[
\norm{f_j} \le \kappa \enpar{ \norm{f_j-f} + \norm{f}}, \quad j\in\NN,
\]
gives the left-hand side inequality. In turn, the inequality
\[
\norm{f} \le \kappa \enpar{ \norm{f-f_j} + \norm{f_j}}, \quad j\in\NN,
\]
gives the right-hand side one.
\end{proof}

The \emph{Banach envelope} of a quasi-Banach space $\XX$ consists of a Banach space $\widehat{\XX}$ together with a linear contraction $J_\XX\colon \XX \to \widehat{\XX}$ satisfying the following universal property: for every Banach space $\YY$ and every continuous linear map $T\colon \XX \to \YY$ there is a unique continuous linear map $\widehat{T}\colon \widehat{X}\to Y$ such that $\widehat{T}\circ J_\XX=T$, that is, the diagram
\[
\xymatrix{
\widehat{\XX}\ar[drr]^{\widehat{T}} & \\
\XX \ar[u]^{J_\XX} \ar[rr]_T && \YY
}
\]
commutes, and the norm of ``extension'' $\widehat{T}$ is bounded by the norm of $T$. In particular, $\XX$ and $\widehat{\XX}$ have the same dual space. We will refer to $J_{\XX}$ as the \emph{Banach envelope map} of $\XX$. Since the Banach envelope of $\XX$ is defined by means of a universal property, it is unique in the sense that if a Banach space $\UU$ and a bounded linear map $J\colon \XX\to \UU$ satisfy the property, then there is an linear isometry $S\colon \widehat{\XX} \to \UU$ with $S\circ J_\XX=J$.

To witness the existence of the Banach envelope of $\XX$, we can appeal to the \emph{bidual map} $h_\XX\colon\XX \to \XX^{**}$. In fact, we can choose $\widehat{\XX}=\overline{h_\XX(\XX)}$ and $J_\XX=h_\XX$. A more constructive, and often more useful, approach to Banach envelopes goes as follows. Let $\norm{\cdot}_c$ be the Minkowski functional of the convex hull of $B_\XX$, and put
\[
N=\enbrace{f\in \XX \colon \norm{f}_c=0}.
\]
The Banach envelope of $\XX$ is the completion of the normed space $(\XX/N, \norm{\cdot}_c)$ (see e.g.\@ \cite{AACD2018}).

We point out that the Banach envelope map $J_\XX$ might not be one-to-one. In fact, $J_\XX$ is one-to-one if and only if $\XX$ has the \emph{separation property}, that is, its dual space $\XX^*$ separates the points of $\XX$.

Let $\delta_{\cdot,\cdot}$ denote the Kronecker delta. A \emph{minimal system} in a quasi-Banach space $\XX$ is a sequence $\XB=(\xx_n)_{n=1}^\infty $ in $\XX$ for which there is $\XB^*=(\xx_n^*)_{n=1}^\infty $ in $\XX^*$ such that
\[
\xx_n^*(x_k)=\delta_{n,k}, \quad n,k\in\NN.
\]
We say that $\XB^*$ is a sequence of \emph{coordinate functionals} for $\XB$. If $(\xx_n^*)_{n=1}^\infty $ is bounded, then the \emph{coefficient transform}
\[
\Ft[\XB]\colon\XX\to \FF^{\NN}, \quad
f\mapsto\enpar{\xx_n^*(f)}_{n=1}^\infty , \quad f\in\XX,
\]
is a bounded map from $\XX$ to $\ell_\infty$.

If the closed linear span of $\XB$ is the whole space $\XX$ then the coordinate functionals are unique and we say that $\XB$ is \emph{complete}. If $\XB$ is complete and $\spn(\XB^*)$ is weak*-dense, we say that $\XB$ is \emph{total} and call it a \emph{Markushevich basis}. Note that a complete minimal system $\XB$ is total if and only if the coefficient transform is one-to-one. Hence, the existence of a Markushevich basis for $\XX$ guarantees that $\XX$ has the separation property.

Given a complete minimal system $\XB=(\xx_n)_{n=1}^\infty $ with coordinate functionals $(\xx_n^*)_{n=1}^\infty $, and $A\in\enbrak{\NN}^{<\infty}$, the \emph{coordinate projection} relative to $\XB$ associated with $A$ is the linear operator
\[
S_A[\XB]\colon\XX\to\XX, \quad f\mapsto \sum_{n\in A} \xx_n^*(f) \, \xx_n.
\]
More generally, given $\tau=(t_n)_{n=1}^\infty\in \cc_{00}$, we consider the \emph{multiplier operator}
\[
M_\tau[\XB]\colon\XX \to \XX, \quad f\mapsto \sum_{n=1}^\infty t_n \, \xx_n^*(f) \, \xx_n.
\]
We say that $\XB$ is \emph{$M$-bounded} if
\[
\sup_{n\in\NN} \norm{S_{\{n\}}[\XB]}=\sup_{n\in\NN} \norm{\xx_n} \norm{\xx_n^*}<\infty,
\]
and we say that $\XB$ is \emph{semi-normalized} if
\[
\inf_{n\in\NN} \norm{\xx_n}>0, \quad \sup_{n\in\NN} \norm{\xx_n} <\infty.
\]
Note that $\XB$ is semi-normalized and $M$-bounded if and only if both $\XB$ and $\XB^*$ are bounded.

An \emph{unconditional basis} of a (separable) quasi-Banach space $\XX$ is a sequence $(\xx_n)_{n=1}^\infty $ in $\XX$ such that for every $f\in\XX$ there is a unique sequence $(a_n)_{n=1}^\infty $ in $\FF$ for which the series $\sum_{n=1}^\infty a_n\, \xx_n$ unconditionally converges fo $f$. It is known (see e.g.\@ \cite{AABW2021}) that $(\xx_n)_{n=1}^\infty $ is an unconditional basis if and only if it is a complete minimal system with
\[
\sup\enbrace{\norm{S_A[\XB]} \colon A\in\enbrak{\NN}^{<\infty}}<\infty.
\]
Also, $(\xx_n)_{n=1}^\infty $ is an unconditional basis if and only if it is a complete minimal system with
\[
C_u:=\sup\enbrace{\norm{M_\tau[\XB]} \colon \tau\in \cc_{00},\, \norm{\tau}_\infty\le 1}<\infty.
\]
If $C_u\le C<\infty$, we say that $\XB$ is $C$-unconditional.

A Schauder basis of a (separable) quasi-Banach space $\XX$ is a sequence $(\xx_n)_{n=1}^\infty $ in $\XX$ such that for every $f\in\XX$ there is a unique scalar sequence $(a_n)_{n=1}^\infty $ for which the series $\sum_{n=1}^\infty a_n\, \xx_n$ converges to $f$. Note that any arrangement of an unconditional basis is a Schauder basis.

We denote by $\It$ the set of all finite integer intervals contained in $\NN$. It is known (see e.g.\@ \cite{AABW2021}) that $(\xx_n)_{n=1}^\infty $ is a Schauder basis of a quasi-Banach space if and only if it is a complete minimal system with
\[
\sup_{I\in\It} \norm{S_I[\XB]}<\infty.
\]
So, Schauder bases are in particular $M$-bounded Markushevich bases.

A minimal system $(\xx_n)_{n=1}^\infty $ is said to be \emph{$C$-symmetric}, $1\le C<\infty$, if
\[
\norm{\sum_{n=1}^\infty a_n \, \xx_n} \le C \norm{\sum_{n=1}^\infty a_n \, \xx_{\pi(n)}}, \quad (a_n)_{n=1}^\infty \in \cc_{00},
\]
for every permutation $\pi$ of $\NN$. We say $\XB$ is symmetric if it is $C$-symmetric for some $C$.

To study greedy-like properties, we restrict ourselves to minimal systems that are complete, semi-normalized, and $M$-bounded. For convenience, we will refer to such systems as \emph{$W$-bases}. Let $\XB=(\xx_n)_{n=1}^\infty $ be a $W$-basis of a (separable) quasi-Banach space $\XX$. A set $A\in \enbrak{\NN}^{<\infty}$ is said to be a \emph{greedy set} of $f\in\XX$ relative to $\XB$ if
\[
\abs{\xx_n^*(f)} \ge \abs{\xx_k^*(f)}, \quad n\in A, \, k\in\NN\setminus A.
\]
We denote by $\Gr[\XB](f)$ the set of all greedy sets of $f$. Note that $\emptyset\in\Gr[\XB](f)$ for all $f\in\XX$. Following \cite{KoTe1999}, we say that $\XB$ is quasi-greedy if
\[
\sup\enbrace{ \norm{f-S_A[\XB](f)} \colon f\in B_\XX, \, A\in\Gr[\XB](f)}<\infty.
\]

Let $\XB$ be a $W$-basis of a quasi-Banach space $\XX$. Pick $f\in\XX$. Since $\Ft[\XB](f)\in \cc_{0}$, the set $\Gr[\XB](f)$ is directed by inclusion, that is, for every $A$, $B\in\Gr[\XB](f)$ there is $D\in\Gr[\XB](f)$ with $A\cup B\subset D$. So, we can consider limits of families modelled over $\Gr[\XB](f)$. Note that, since $\Gr[\XB](f)\cap\enbrak{A}$ is finite for all $A\in\Gr[\XB](f)$, convergent nets modelled over $\Gr[\XB](f)$ are bounded. Also note that, given a directed set $\Lambda$, a monotone function $h\colon \Lambda \to \Gr[\XB](f)$ is final if and only if $\lim_{\lambda\in\Lambda} \abs{h(\lambda)}=\infty$.

In this terminology, the characterization of quasi-greedy bases as those bases for which the TGA converges (see \cites{Woj2000,AABW2021}) can be restated as follows.

\begin{theorem}[cf.\@ \cite{AABW2021}*{Theorem 4.1}]\label{thm:AABWCha}
A $W$-basis $\XB$ of a quasi-Banach space $\XX$ is quasi-greedy if and only
\[
\lim_{A\in\Gr[\XB](f)} S_A[\XB](f)=f
\]
for all $f\in\XX$.
\end{theorem}

We point out that Theorem~\ref{thm:AABWCha} implies that any quasi-greedy basis is total (see \cite{AABW2021}*{Corollary 4.5}).

Following \cite{DKKT2003}, we say that a $W$-basis is \emph{almost greedy} if there is a constant $C\in[1,\infty)$ such that
\begin{equation*}
\norm{f-S_A[\XB](f)} \le C \norm{f-S_B[\XB](f)}, \quad A\in\Gr[\XB](f), \, B\in\enbrak{\NN}^{\le\abs{A}}.
\end{equation*}

In greedy approximation, several nonlinear forms of unconditionality and symmetry naturally appear. To properly define them, we introduce some terminology. Let $\XB=(\xx_n)_{n=1}^\infty $ be a sequence in $\XX$. Given $A\in\enbrak{\NN}^{<\infty}$ finite and $\varepsilon=(\varepsilon_n)_{n\in A}\in S_\FF^A$, we consider the \emph{indicator vector} and the \emph{signed indicator vector}
\[
\Ind_{A}[\XB]=\sum_{n\in A} \xx_n, \quad
\Ind_{\varepsilon,A}[\XB]=\sum_{n\in A} \varepsilon_n \, \xx_n.
\]
We say that $\XB$ is \emph{unconditional for constant coefficients} if there is a constant $C\in[1,\infty)$ such that
\begin{equation*}
\norm{\Ind_{\varepsilon,A}[\XB]} \le C \norm{\Ind_{\varepsilon,B}[\XB]}, \quad B\in\enbrak{\NN}^{<\infty}, \, A\subset B, \, \varepsilon\in S_\FF^B.
\end{equation*}
Unconditional bases are quasi-greedy. In turn, quasi-greedy bases are unconditional for constant coefficients \cite{Woj2000}.

We say that $\XB$ is \emph{super-democratic} if there is a constant $C$ such that
\begin{equation*}
\norm{\Ind_{\delta,A}[\XB]} \le C \norm{\Ind_{\varepsilon,B}[\XB]}, \quad
B\in\enbrak{\NN}^{<\infty}, \, A\in\enbrak{\NN}^{\le \abs{B}}
\, \delta\in\EE^A,\, \varepsilon\in\EE^B,
\end{equation*}
and we say that $\XB$ is \emph{democratic} if there is a constant $C$ such that
\begin{equation*}
\norm{\Ind_{A}[\XB]} \le C \norm{\Ind_{B}[\XB]}, \quad B\in\enbrak{\NN}^{<\infty}, \, A\in\enbrak{\NN}^{\le \abs{B}}.
\end{equation*}
Simmetry implies democracy. Superdemocracy can be characterized in terms of the \emph{upper democracy function}, also known as \emph{fundamental function}, $\usdf$ and the \emph{lower democracy function} $\lsdf$ of the basis. Namely, if for $m\in\NN$ we set
\begin{align*}
\usdf(m) &=\usdf[\XB](m)=\sup \enbrace{ \norm{\Ind_{\varepsilon, A}[\XB]}, \, A\in\enbrak{\NN}^{\le m}, \,\varepsilon\in\EE^A} ,\\
\lsdf(m) &=\lsdf[\XB](m) =\inf \enbrace{ \norm{\Ind_{\varepsilon, A}[\XB]}, \, A\in\enbrak{\NN}^{\ge m}, \,\varepsilon\in\EE^A},
\end{align*}
then $\XB$ is superdemocratic if and only $\usdf(m) \lesssim \lsdf(m)$ for $m\in\NN$. Also, $\XB$ is superdemocratic if and only if it is democratic and unconditional for constant coefficients. Note that $\XB$ is semi-normalized if and only is $\lsdf(1)>0$ and $\usdf(1)<\infty$, in which case $\lsdf(m)>0$ and $\usdf(m)<\infty$ for all $m\in\NN$.

While symmetric Schauder bases are unconditional, there are Banach spaces with conditional symmetric minimal systems (see \cite{AABCO2024}*{Example 3.9}). In fact, we have the following.

\begin{theorem}[\cite{AABCO2024}*{Corollary 3.8}]\label{thm:symunc}
Given a symmetric complete minimal system $\XB$ of a Banach space $\XX$, the following are equivalent.
\begin{enumerate}[label=(\roman*), leftmargin=*]
\item $\XB$ is unconditional.
\item $\XB$ is unconditional for constant coefficients.
\item $\XB$ is total.
\end{enumerate}
\end{theorem}

For further reference, we record the well-known characterization of almost greedy bases.

\begin{theorem}[\cite{AABW2021}*{Theorem 6.3}]\label{thm:CAG}
A $W$-basis is almost greedy if and only if it is quasi-greedy and democratic.
\end{theorem}

We point out that Theorem~\ref{thm:CAG} was originally stated and proved in \cite{DKKT2003} for Schauder bases in Banach spaces. Generalizing it to W-bases in quasi-Banach spaces requires new techniques.

If $\XB=(\xx_n)_{n=1}^\infty $ is a minimal system of $\XX$, then
\[
\widehat{\XB}:=\enpar{J_\XX(\xx_n)}_{n=1}^\infty
\]
is a minimal system of the Banach envelope $\widehat{\XX}$ of $\XX$. Those properties of $\XB$ that can be characterized in terms of boundedness of linear operators plainly pass from $\XB$ to its Banach envelope $\widehat{\XB}$. For instance, if $\XB$ is either Schauder, unconditional, symmetric, complete, or $M$-bounded, so is $\widehat{\XB}$. We also note that if $\XB$ is semi-normalized and $M$-bounded, then $\widehat{\XB}$ is semi-normalized. In contrast, as we will see, there are instances where $\XB$ is total but $\widehat{\XB}$ is not.

Let $(\ee_n)_{n=1}^\infty $ be the sequence of \emph{unit vectors}, that is, for each $n\in\NN$
\[
\ee_n\colon\NN\to\{0,1\}
\]
is given by $\ee_n(k)=\delta_{n,k}$ for all $k\in\NN$. The \emph{unit functionals} $\ee_n^*$, $n\in\NN$, will be the scalar-valued linear maps on $\FF^\NN$ given by
\[
(a_k)_{k=1}^\infty\mapsto a_n.
\]

Some of the notions we have introduced related to minimal systems make sense for the unit vector system $\EB=(\ee_n)_{n=1}^\infty $ regardless of the existence of a topology on $\cc_{00}$. The support of $f\in\FF^\NN$ is the set
\[
\supp(f)=\enbrace{n\in\NN \colon \ee_n^*(f)\not=0}.
\]
The coordinate projection $S_A\colon\FF^\NN \to \FF^\NN$ is well-defined for all $A\in\enbrak{\NN}$, and the multiplier $M_\tau\colon\FF^\NN \to \FF^\NN$ is well-defined for all $\tau\in\FF^{\NN}$. Given a one-to-one map $\pi\colon\NN\to\NN$ we define
\[
P_\pi\colon \FF^{\NN} \to \FF^{\NN}, \quad f=(a_n)_{n=1}^\infty \mapsto (b_k)_{k=1}^\infty, \quad
b_k=\begin{cases}
a_n & \mbox{ if } k=\pi(n),\\ 0 & \mbox{ if }k\in\NN\setminus\pi(\NN).
\end{cases}
\]
Note that $P_\pi$ is a bijection if and only if $\pi$ is a permutation of $\NN$. We will denote by $\Ind_A$ and $\Ind_{\varepsilon,A}$, $A\in\enbrak{\NN}^{<\infty}$, $\varepsilon\in S_\FF^A$, the corresponding indicator and signed indicator vectors relative to $\EB$, and by $\Ind_A^*$ and $\Ind_{\varepsilon,A}^*$ those relative to $(\ee_n^*)_{n=1}^\infty $. We say that $A\in\enbrak{\NN}^{<\infty}$ is a greedy set of $f$ if $\abs{a_n}\le \abs{a_k}$ for all $n\in \NN\setminus A$ and $k\in A$. We denote by $\Gr(f)$ the set consisting of all greedy sets of $f$. While there are sequences whose unique greedy set is the empty set, $\Gr(f)$ is directed by inclusion for every $f\in\cc_0$, and for each $m\in\NN$ there exists $A\in\Gr(f)$ with $\abs{A}=m$.

Given $t\in\FF$ we define $\sgn(t)\in S_\FF$ by $\sgn(t)=t/\abs{t}$ if $t\not=0$ and $\sgn(t)=1$ otherwise. Given $f=(a_n)_{n=1}^\infty \in\FF^{\NN}$ we put
\[
\varepsilon(f)=\enpar{\sgn(a_n)}_{n=1}^\infty .
\]

A \emph{sequence space} will be a quasi-Banach space $\Sym$ such that $\cc_{00}\subset \Sym$, $\norm{\ee_n}_\Sym\approx 1$ for $n\in\NN$, and $\Sym\subset \ell_\infty$ continuously. If $\cc_{00}$ is dense in $\Sym$, in which case $\Sym\subset \cc_0$, we say that $\Sym$ is \emph{minimal}. In other words, a minimal sequence space is a quasi-Banach space for which the unit vectors form a semi-normalized $M$-bounded Markushevich basis.

Given a sequence space $\Sym$ we denote by $\lsdf[\Sym]$ (resp.\@ $\usdf[\Sym]$) the lower (resp.\@ upper) democracy function of its unit vector system. Given a linear operator $T\colon\FF^{\NN} \to \FF^{\NN}$ and sequence spaces $\Sym_j$, $j=1$, $2$,
\[
\norm{T}_{\Sym_1\to\Sym_2}
\]
will denote the norm of $T$ regarded as an operator from $\Sym_1$ into $\Sym_2$.

A sequence space $\Sym$ is said to be unconditional if $\norm{M_\tau}_{\Sym\to\Sym}\le \norm{\tau}_\infty$ for every $\tau\in\ell_\infty$. In other words, an unconditional sequence space is a quasi-Banach lattice modelled over $\NN$ such that the unit vector system is semi-normalized. If $\norm{P_\pi}_{\Sym\to\Sym}\le 1$ for every permutation $\pi$ of $\NN$, we say that the sequence space $\Sym$ is symmetric. Clearly, the unit vector system of an unconditional (resp., symmetric) sequence space is a $1$-unconditional (resp., $1$-symmetric) basis of the separable part of the space.

Given $f\in \cc_{00}$, $f\ge 0$, we denote by $D(f)$ its nonincreasing rearrangement. Then, we extend $D$ to an operator from $\FF^\NN$ to $[0,\infty]^\NN$ by
\[
D(f)=\sup \enbrace{D(g) \colon 0 \le g\le \abs{f}, \, g\in \cc_{00}}.
\]
Given $1<q\le \infty$, $k\in\uii{0}$ and $f\in\FF^\NN$ we define
\[
\rho_{1,q}(f,k)=\norm{ \enpar{m^{1-1/q} D(f)(m+k)}_{m=1}^\infty}_q.
\]
Let $\ell_{1,q}$ denote the space consisting of all $f\in\FF^\NN$ such that
\[
\norm{f}_{1,q}:=\rho_{1,q}(f,0)<\infty.
\]

The Lorentz sequence space $\ell_{1,q}$ is a symmetric unconditional sequence space with modulus of concavity at most $2$ (see e.g.\@ (\cite{CRS2007}*{Proof of Lemma 2.2.10}). Besides, $\ell_{1,q}$, while fails to be locally convex, is locally $p$-convex for any $0<q<1$ (see \cite{Hunt1966}*{Statements~(2.2) and (2.6)}). Concerning its fundamental function, we have
\[
\norm{\Ind_{\varepsilon,A}}_{1,q}\approx\abs{A}, \quad A\in\enbrak{\NN}^{<\infty}, \, \varepsilon\in S_\FF^A.
\]
The sequence space $\ell_{1,q}$ is minimal if only if $q<\infty$ (see \cite{Hunt1966}*{Statement~(2.4)}). The separable part of $\ell_{1,\infty}$, denoted by $h_{1,\infty}$, can be expressed as
\[
\enbrace{ f\in\FF^\NN \colon \lim_{k\in\NN} \rho_{1,\infty}(f,k)=0}
=\enbrace{ f\in\FF^\NN \colon \lim_{m\in\NN} m D(f)(m)=0}.
\]
To deal with the space $h_{1,\infty}$, we will use an inequality involving the map $\rho_{1,\infty}$ that relies on the well-known estimate
\begin{equation}\label{eq:BS}
D(f+g)(m+n) \le D(f)(m)+D(g)(n), \quad f,g\in\FF^\NN, \, m,n\in\NN
\end{equation}
(see \cite{Hunt1966}*{Statement~(1.6)}).

\begin{lemma}\label{lem:BS}
For any $f$, $g\in\FF^\NN$, and $k_1$, $k_2\in\uii{0}$,
\[
\rho_{1,\infty}(f+g,1+k_1+k_2)\le 2 \enpar{ \rho_{1,\infty}\enpar{f,k_1}+\rho_{1,\infty}\enpar{g,k_2}}.
\]
\end{lemma}

\begin{proof}
Fix $m\in\NN$ and set $j=\enfloor{(m+1)/2}$. Since $m\ge 2j-1$, by \eqref{eq:BS},
\begin{align*}
m D(f+g)(m+k_1+k_2+1)&\le m D(f+g)(2j+k_1+k_2)\\
&\le m\enpar{D(f)(j+k_1) +G(g)\enpar{j+k_2}}\\
&\le \frac{m}{j} \enpar{ \rho_{1,\infty}\enpar{f,k_1}+\rho_{1,\infty}\enpar{g,k_2}}.
\end{align*}
Since $m\le 2j$, we are done.
\end{proof}

We will also use the following embedding result from \cite{AABW2021}.

\begin{theorem}[see \cite{AABW2021}*{Section 9}]\label{thm:LorentzEmbeds}
Let $\Sym$ be a sequence space whose canonical basis is quasi-greedy. If $m\lesssim \lsdf[\Sym]$ for $m\in\NN$, then $\Sym\subset h_{1,\infty}$ continuously.
\end{theorem}

We record for further reference three elementary lemmas.

\begin{lemma}\label{lem:GUnion}
Let $f\in\FF^\NN$ and $A$, $B\in\Gr(f)$. Then $A\cup B\in\Gr(f)$.
\end{lemma}

\begin{proof}
It is clear from the definition.
\end{proof}

\begin{lemma}\label{lem:GSProj}
Let $f\in\FF^\NN$ and $A\in\Gr(f)$ and $B\in\Gr(f-S_A(f))$. Then $A\cup B\in\Gr(f)$.
\end{lemma}

\begin{proof}
Pick $n\in A\cup B$ and $k\in\NN\setminus(A\cup B)$. If $n\in A$ then, since $k\in\NN\setminus A$, $\abs{\ee_k^*(f)}\le \abs{\ee_n^*(f)}$. If $n\notin A$, then since $k\in\NN\setminus B$,
\[
\abs{\ee_k^*(f)}=\abs{\ee_k^*(f-S_A(f))}\le \abs{\ee_n^*(f-S_A(f))} =\abs{\ee_n^*(f)}.\qedhere
\]
\end{proof}

\begin{lemma}\label{lem:B99}
Given $f\in\FF^\NN$, $A\in\Gr(f)$, $B\in\enbrak{\NN\setminus A}^{<\infty}$ and $k\in\bii{0}{\abs{A}}$,
\[
\abs{\Ind^{*}_B(f)}\le \frac{\abs{B}}{1+\abs{A}-k} \rho_{1,\infty}(f,k)
\]
\end{lemma}

\begin{proof}
Set $m=\abs{A}$. Since $\abs{\ee_n^*(f)}\le D(f)(m+1)$ for all $n\in B$,
\[
\abs{\Ind^{*}_B(f)}
\le \abs{B} D(f)(m+1)
\le \frac{\abs{B}}{1+m-k} \rho_{1,\infty}(f,k).\qedhere
\]
\end{proof}
%--------------------------------------------------------------------------
\section{Quasi-greediness does not pass to envelopes}\label{sect:AGDNPE}\noindent
%--------------------------------------------------------------------------
Our construction relies on the introduction of a remarkable space related to the TGA. We define $\norm{\cdot}_{\BB}\colon\FF^\NN\to[0,\infty]$ by
\[
\norm{f}_\BB=\sup\enbrace{ \abs{\Ind^{*}_{I\setminus A} (f)} \colon A\in\Gr(f), \, I\in\It}.
\]

Since $\Gr(tf)=\Gr(f)$ for all $f\in\FF^\NN$ and $t\in\FF\setminus\{0\}$, $\norm{t\,f}_\BB=\abs{t} \norm{f}_\BB$ for all $t\in\FF$ and $f\in\FF^{\NN}$. In contrast, as we will see,
\[
\BB=\enbrace{f\in\FF^\NN \colon \norm{f}_\BB<\infty}
\]
is not a vector space. We define for each $f\in\FF^\NN$ and $A\in\enbrak{\NN}$
\[
\beta(f,A)=\sup_{I\in\It} \abs{\Ind^{*}_{I\setminus A} (f)},
\]
and we consider the subspace of $\BB$ defined as
\[
\BB_0=\enbrace{f\in \cc_0 \colon \lim_{A\in\Gr(f)} \beta(f,A)=0}.
\]

We start our study of $\BB$ by estimating $\norm{\cdot}_\BB$ in some instances.

\begin{lemma}\label{lem:Schauder}
For all $f\in\FF^\NN$ and all integer intervals $J$ we have
$\norm{S_J(f)}_\BB\le \norm{f}_\BB$.
\end{lemma}

\begin{proof}
Just note that $J\cap I\in\It$ for every $I\in\It$.
\end{proof}

\begin{lemma}\label{lem:B0}
For all $f\in\FF^\NN$ we have $\norm{f}_\BB\le\norm{f}_1$. If $f$ is nonnegative, then $\norm{f}_\BB=\norm{f}_1$. Besides, $\ell_1\subset\BB_0$.
\end{lemma}

\begin{proof}
Given $A\in\Gr(f)$ and $I\in\It$,
\[
\abs{\Ind^{*}_{I\setminus A}(f)}
\le \Ind^{*}_{ I \setminus A}(\abs{f})
\le \norm{f-S_A(f)}_1.
\]
Hence, $\beta(f,A)\le \norm{f-S_A(f)}_1$. By unconditionality and Theorem~\ref{thm:AABWCha}, $\norm{f}_\BB\le \norm{f}_1$, and $f\in\BB_0$ as long as $\norm{f}_1<\infty$. If $f\ge 0$, since $\emptyset\in\Gr(f)$,
\[
\norm{f}_\BB\ge\sup_{I\in\It} \Ind^{*}_I(f)=\norm{f}_1.\qedhere
\]
\end{proof}

A sequence $(J_k)_{k=1}^\infty$ in $\enbrak{\NN}$ is said to be \emph{right-dominant} if $J_k\not=\emptyset$ and
\[
\max(J_k)<\min(J_{k+1})
\]
for all $k\in\NN$. A nonnegative sequence $g=(a_n)_{n=1}^\infty$ is said to be \emph{Leibnizian} if there is a right-dominant sequence $\Jt:=(J_k)_{k=1}^\infty$ in $\It$ such that
\begin{itemize}
\item $\supp(g)\subset\cup_{k=1}^\infty J_k$,
\item $(\Ind^*_{J_k}(g))_{k=1}^\infty$ is nonincreasing, and
\item $\min\enbrace{a_n \colon n\in J_k\cap\supp(g)}>\max\enbrace{a_n \colon n\in J_{k+1}}$.
\end{itemize}
Such a sequence $\Jt$ is said to be \emph{admissible} for $g$. Using the convention that $\min(\emptyset)=\infty$, the existence of $k\in\NN$ such that $J_k\cap\supp(g)=\emptyset$ is not excluded. If it is the case and the set
\[
S(g,\Jt):=\enbrace{k\in\NN \colon J_k\cap\supp(g)\not=\emptyset}=\enbrace{k\in\NN \colon \Ind^*_{J_k}(g)\not=0}
\]
is nonempty, then there is $m\in\NN$ such that $S(g,\Jt)=\bii{1}{m}$. Any $g\in\cc_{00}$ is Leibnizian.

If $g$ is Leibnizian and $\Jt:=(J_k)_{k=1}^\infty$ is admissible for it, we set
\[
\alpha(g,\Jt)=\Ind^*_{J_1}(g), \quad \omega(g, \Jt)=\lim_k \Ind^*_{J_k}(g).
\]

\begin{lemma}\label{lem:LeibB}
If $g\in[0,\infty)^\NN$ is decreasing on its support then $g$ is Leibnizian. In fact, any right-dominant sequence $\Jt=(J_k)_{k=1}^\infty$ in $\It$ such that $J_k\cap\supp(g)$ is either a singleton or empty for all $k\in\NN$ is admissible for $g$, and $\alpha(g,\Jt)=\norm{g}_\infty$.
\end{lemma}

\begin{proof}
It is clear from definition.
\end{proof}

\begin{lemma}\label{lem:LeibA}
Let $t>0$ and $g\in \cc_0 \setminus \ell_1$ be nonnegative. Then there is $A\in\enbrak{\NN}$ such that $S_A(g)$ is Leibnizian, and there is $\Jt$ admissible for $g$ such that $\omega(g,\Jt)\ge t$.
\end{lemma}

\begin{proof}
Pick closed intervals $(H_k)_{k=1}^\infty$ contained in $[t,\infty)$ such that
\[
\max(H_{k+1})<\min(H_k)<\max(H_k)
\]
for all $k\in\NN$. We recursively construct a right-dominant sequence $(J_k)_{k=1}^\infty$ in $\It$ such that, with the convention that $J_0=\emptyset$,
\[
\min\enbrace{a_n \colon n\in J_{k-1}\cap\supp(g)}>\max\enbrace{a_n \colon n\in J_k},
\]
and $\Ind_{J_k}^*(g)\in H_k$ for all $k\in\NN$. The set $A=\cup_{k=1} J_k$ satisfies the desired conditions.
\end{proof}

\begin{lemma}\label{lem:B77}
Let $g$ be a Leibnizian sequence and $\Jt=(J_k)_{k=1}^\infty$ be admissible for $g$.
\begin{enumerate}[label=(\roman*),leftmargin=*,widest=ii]
\item Let $\tau=(t_n)_{n=1}^\infty\in\{-1,1\}^\NN$ be such that $t_n=(-1)^{k-1}$ for all $k\in\NN$ and $n\in J_k$. Put $f=M_\tau(g)$. Then $\norm{f}_\BB=\alpha(g,\Jt)$ and
\[
\lim_{A\in\Gr(f)} \beta(f,A)=\omega(g,\Jt).
\]
In particular, if $\omega(g,\Jt)=0$ then $f\in\BB_0$.

\item If $\omega(g,\Jt)>0$, then there is an uncountable family $(\rho_j)_{j\in J}$ in the unit ball of $\ell_\infty$ such that $\norm{M_{\rho_j}(g)}_\BB\le \alpha(g,\Jt)$ for all $j\in J$, and
\[
\norm{M_{\rho_j}(g) -M_{\rho_k}(g) }_\BB\ge\omega(g,\Jt), \quad j,\, k\in J, \, j\not=k.
\]
\end{enumerate}
\end{lemma}
\begin{proof}
Clearly, $A\in\Gr(f)$ if and only if there are $k=k(A)\in\NN$ and $B=B(A) \in\Gr(S_{J_k}(g))$
such that
\[
A=\enpar{\supp(g) \cap \enpar{\cup_{j=1}^{k-1} J_j}} \cup B.
\]
Thus, given $I\in\It$, there are $j$, $m\in\NN$, $D\in\enbrak{J_j\cap \supp(g)}$ and $E\in\enbrak{J_{m+1}\cap \supp(g)}$ such that $k(A)\le j\le m$ and
\[
\supp(f)\cap I\setminus A=D\cup \enpar{\cup_{k=j+1}^m J_k \cap \supp(g)} \cup E.
\]
Hence, there are $\bii{j+1}{m}\subset J\subset\bii{j}{m+1}$ and $\varepsilon\in\{-1,1\}$ so that
\[
\abs{\Ind_{I\setminus A}^*(f)}\le \varepsilon\sum_{k\in J} (-1)^k \Ind_{J_k}^*(g).
\]
By Leibniz' criterion, $\abs{\Ind_{I\setminus A}^*(f)}\le \Ind_{J_j}^*(g)$. Therefore,
\[
\beta(f,A)\le \Ind_{J_{k(A)}}^*(g).
\]
Consequently, $\norm{f}_\BB\le \alpha(g,\Jt)$ and, since $\lim_{A\in\Gr(f)} k(A)=\infty$,
\[
\limsup_{A\in\Gr(f)} \beta(f,A)\le \omega(g,\Jt).
\]
Given $k\in\NN$, let $A_k\in\Gr(f)$ by such that $k(A)=k$ and $B(A)=\emptyset$, and set $I_k=\bii{1}{\max(J_k)}$. We have
\[
\beta(f,A_k)\ge \abs{\Ind^*_{I_k\setminus A_k}(f)}=\Ind_{J_k}^*(g).
\]
Therefore, $\norm{f}_\BB\ge \alpha(g,\Jt)$ and, since $(A_k)_{k=1}^\infty$ is a monotone final sequence,
\[
\liminf_{A\in\Gr(f)} \beta(f,A)\ge \omega(g,\Jt).
\]

Given $N\in\enbrak{\NN}$, set $A_N=\cup_{k\in N} J_k$. The sequence
\[
g_N= S_{A_N}(g)
\]
is Leibnizian, and there is $\Jt_N:=(J_{k,N})_{k=1}^\infty$ admissible for $g_N$ so that for every $k\in S(g_N,\Jt_N)$ there is $j\in N$ such that $J_{k,N}=J_j$. If $f_N$ and $\tau_N$ are as before relative to $g_N$ and $\Jt_N$, then $f_N=M_{\rho_N}(g)$, where
\[
\rho_N=\tau_N \chi_{A_N}.
\]
If $N\not=\emptyset$, then
\begin{align}
\Ind^*_{J_{1,N}}(g_N)=\Ind^*_{J_{\min(N)}}(g).\label{m1}
\end{align}
Hence, $\norm{f_N}_{\BB}\le \alpha(g,\Jt)$ for all $N\in\enbrak{\NN}$. To complete the proof, we note that for each $N\in\enbrak{\NN}$ there is $(\varepsilon_{k,N})_{k=1}^\infty$ such that $\varepsilon_{k,N}\in\{-1,1\}$ if $k\in N$, $\varepsilon_{k,N}=0$ if $k\notin N$, and
\[
S_{J_k}(f_N)=\varepsilon_{k,N} S_{J_k}(g), \quad k\in\NN.
\]
Therefore, given $N$, $M\in\enbrak{\NN}$ with $N\not=M$,
\[
\norm{f_N-f_M}_{\BB} \ge \sup_{k\in M\triangle N} \abs{\Ind^*_{J_k}(f_N-f_M)}=\sup_{k\in M\triangle N} \Ind^*_{J_k}(g)\ge \omega(g,\Jt).\qedhere
\]
\end{proof}

\begin{lemma}\label{lem:BNorm3}
Let $f=(a_n)_{n=1}^\infty $ be a sequence in $[0,1]$. Assume that $a_1=a_2=0$, and $a:=\sum_{n=1}^\infty a_n\in[1,\infty)$. Then,
\[
\norm{\ee_1- t \ee_2+f}_\BB
=\begin{cases}
1+a & \mbox{ if } t= 1, \\
1+a -t& \mbox{ if } 0\le t<1.
\end{cases}
\]
\end{lemma}

\begin{proof}
Set $g=\ee_1- t \ee_2+f$ and $h=\ee_1+f$. On one hand,
\[
-t \le \Ind^{*}_B(g)\le 1+a
\]
for all $B\in\enbrak{\NN}^{<\infty}$. On the other hand,
\[
\sup\enbrace{\Ind^{*}_I(g) \colon I\in\It, \min(I)=1}=1+a-t.
\]
Hence,
\[
a+1-t\le \norm{g}_\BB \le \max\enbrace{t,1+a}=	1+a.
\]

If $t=1$, $\{2\}$ is a greedy set of $f$. Therefore,
\[
\norm{g}_\BB\ge \sup\enbrace{ \abs{\Ind^{*}_{I\setminus\{2\}} (g)} \colon I \in\It} =\sup\enbrace{ \Ind^{*}_I(h) \colon I \in\It}=1+a.
\]
Pick $0<t<1$, $A\in\Gr(g)$ and $I\in\It$. If $2\notin B:=I\setminus A$, then $B\subset\{1\}$ or $B\subset[3,\infty)$. Then,
\[
\norm{g}_\BB\le \sup\{1,a,1+a-t\}=1+a-t.\qedhere
\]
\end{proof}

The following lemmas are aimed to give conditions on sequence spaces $\Sym$ which ensure that $\BB\cap \Sym$, equipped with the gauge
\[
f\mapsto \norm{f}_{\BB,\Sym}:=\max\enbrace{ \norm{f}_\BB,\norm{f}_\Sym},
\]
is a quasi-Banach space. For proving the first of them, we will use the identity
\begin{equation}\label{eq:twisted}
\Ind^*_A(f+g)-\Ind^*_{B}(f)-\Ind^*_D(g)
=\Ind^*_{(A\cup D)\setminus B}(f) + \Ind^*_{(A\cup B)\setminus D}(g) -\Ind^*_{(B\cup D)\setminus A}(f+g),
\end{equation}
valid for any $f$, $g\in\FF^\NN$ and $A$, $B$, $D\in\enbrak{\NN}^{<\infty}$, which will be used several more times throughout the paper.

\begin{lemma}\label{lem:B1}
Given $f$, $g\in\FF^\NN$,
\[
\norm{f+g}_\BB\le 2 \enpar{ \norm{f+g}_{1,\infty}+\norm{f}_{1,\infty}+\norm{g}_{1,\infty}}+\norm{f}_\BB+\norm{g}_\BB.
\]
\end{lemma}

\begin{proof}
Pick $A\in\Gr(f+g)$ and $I\in\It$. Choose $B\in\Gr(f)$ and $D\in\Gr(D)$ with $m:=\abs{A}=\abs{B}=\abs{D}$. By \eqref{eq:twisted}, $\Ind^{*}_{I\setminus A}(f+g)$ can be expanded as
\[
\Ind^{*}_{I\setminus B}(f)
+\Ind^{*}_{I\setminus D}(g)
+ \Ind^{*}_{(B\cup D) \cap I\setminus A}(f+g)
- \Ind^{*}_{(A\cup D) \cap I \setminus B}(f)
- \Ind^{*}_{(A\cup B) \cap I \setminus D}(g).
\]
Since $\max\{\abs{B\cup D}, \abs{A\cup D}, \abs{A\cup B}\}\le 2m$, by Lemma~\ref{lem:B99},
\[
\abs{\Ind^{*}_{I\setminus A}(f+g)}
\le \norm{f}_\BB+\norm{g}_\BB+ \frac{2m}{m+1} \enpar{\norm{f}_{1,\infty}+ \norm{g}_{1,\infty} + \norm{f+g}_{1,\infty}}.
\]
Since $2m/(m+1)\le 2$, we are done.
\end{proof}

\begin{lemma}\label{lem:B2}
Suppose $f$ and $(f_j)_{j=1}^\infty$ are in $\FF^\NN$ such that $\lim_j f_j=f$ pointwise. Then
\[
\norm{f}_\BB\le L:= \liminf_{j\in\NN} 6 \norm{f_j}_{1,\infty} + \norm{f_j}_\BB.
\]
\end{lemma}
\begin{proof}
If $\norm{f_j}_{1,\infty}=\infty$ eventually, then $L=\infty$, and we are done. Otherwise, passing to a subsequence we can assume that $\norm{f_j}_{1,\infty}<\infty$, whence $f_j\in\cc_0$, for all $j\in\NN$.

Pick $A\in\Gr(f)$ and $I\in\It$. Assume that $A\setminus I\not=\emptyset$, and set $g=\Ind_{\varepsilon(f),A}$ and $n_0=\max(A\cup I)$. If $t>0$ and $j\in\NN$ is large enough, then there is $B_j\subset [n_0+1, \infty)$ such that $A\cup B_j\in \Gr(f_j+t g)$. Therefore, by Lemma~\ref{lem:B1},
\begin{align*}
&\abs{\Ind^{*}_{I\setminus A}\enpar{f+ t g)}}
=\lim_{j\in\NN} \abs{\Ind^{*}_{I\setminus A}\enpar{f_j + t g}}=\lim_{j\in\NN} \abs{\Ind^{*}_{I\setminus (A\cup B_j)}\enpar{f_j + t g}}\\
&\le \liminf_{j\in\NN} \norm{f_j+ t g}_\BB\\
&\le\liminf_{j\in\NN} 2\enpar{ \norm{f_j+t g}_{1,\infty} + \norm{f_j}_{1,\infty}+ t \norm{g}_{1,\infty}}
+ \norm{f_j}_\BB + t \norm{\Ind_A}_\BB\\
&\le \liminf_{j\in\NN} 6 \enpar{ \norm{f_j}_{1,\infty} + t \norm{g}_{1,\infty}} + \norm{f_j}_\BB + t \norm{g}_\BB.
\end{align*}
By Lemma~\ref{lem:B0}, $\norm{g}_\BB<\infty$. Letting $t$ tend to zero we obtain the desired inequality.
\end{proof}

\begin{remark}\label{rmk:NT}
Suppose that a minimal sequence space $\Sym$ embeds continuously into $\ell_1$. Then, by Lemma~\ref{lem:B0}, $\BB\cap\Sym=\Sym$. Besides, the envelope map of the inclusion map $J\colon \Sym \to \ell_1$ is an isomorphism. Indeed, if we set $\xx_n=\widehat{J}(\ee_n)$ for all $n\in\NN$, there is a bounded linear map $T\colon\ell_1\to \widehat{\Sym}$ such that $T(\ee_n)=\xx_n$ for all $n\in\NN$. We infer that $\widehat{J}$ and $T$ are inverse isomorphisms of one another.
\end{remark}

In light of Remark~\ref{rmk:NT}, we focus on studying $\BB\cap\Sym$ in the case when $\Sym$ is not contained in $\ell_1$. In this regard, we notice that
\[
\ell_1\subsetneq \ell_{1,q} \subsetneq \ell_{1,r}, \quad 1<q<r\le\infty.
\]

\begin{proposition}\label{prop:BSSS}
Let $(\Sym,\norm{\cdot}_\Sym)$ be a minimal sequence space whose canonical basis is quasi-greedy. If $m \lesssim \lsdf[\Sym](m)$ for $m\in\NN$, then $(\BB\cap\Sym,\norm{\cdot}_{\BB,\Sym})$ is a sequence space whose canonical basis is quasi-greedy. Furthermore,
\begin{enumerate}[label=(\roman*),leftmargin=*,widest=iii]
\item\label{it:CCE} $\norm{\Ind_{\varepsilon, A}}_{\BB,\Sym} \approx \norm{\Ind_{\varepsilon, A}}_\Sym$ for $A\in\enbrak{\NN}^{<\infty}$ and $\varepsilon\in S_\FF^A$;
\item\label{it:BSDual} $\norm{\Ind^*_I}_{\BB\cap\Sym\to \FF}\le 1$ for every $I\in\It$;
\item\label{it:PSP} $\norm{S_J(f)}_{\BB} \le \norm{f}_{\BB,\Sym}$ for all $J\in\It$ and $f\in\FF^\NN$; and
\item if $\Sym$ is locally $p$-convex, $0<p\le 1$, then $\BB\cap\Sym$ is locally $q$-convex for all $0<q<p$.
\item If $\Sym$ is unconditional and is not contained in $\ell_1$, then $\BB\cap\Sym$ is nonseparable.
\end{enumerate}
\end{proposition}

\begin{proof}
Assume that $\Sym$ is equipped with a $p$-norm. Set $\kappa=2^{1/p-1}$. Use Theorem~\ref{thm:LorentzEmbeds} to pick $C$ such that $\norm{f}_{1,\infty} \le C \norm{f}_\Sym$ for all $f\in\Sym$. Consider for each $\epsilon>0$ the gauge
\[
\norm{f}_{\BB,\Sym,\epsilon}=\norm{f}_\Sym+\epsilon \norm{f}_\BB, \quad f\in\FF^\NN.
\]
Clearly, $\norm{\cdot}_{\BB,\Sym,\epsilon}$ is equivalent to $\norm{\cdot}_{\BB,\Sym}$. Given $f$, $g\in\FF^\NN$, by Lemma~\ref{lem:B1},
\begin{align*}
&\norm{f+g}_{\BB,\Sym,\epsilon}\\
&\le \enpar{1+ 2C \epsilon} \norm{f+g}_{\Sym} +2C \epsilon \enpar{ \norm{f}_\Sym+\norm{g}_\Sym} + \epsilon\enpar{ \norm{f}_\BB+\norm{g}_\BB}\\
&\le \enpar{\kappa+ (1+\kappa)2C \epsilon} \enpar{ \norm{f}_\Sym+\norm{g}_\Sym} + \epsilon\enpar{ \norm{f}_\BB+\norm{g}_\BB}\\
&\le \kappa_\epsilon \enpar{ \norm{f}_{\BB,\Sym,\epsilon}+\norm{g}_{\BB,\Sym,\epsilon}},
\end{align*}
where $\kappa_\epsilon=\max\{\epsilon,\kappa+ (1+\kappa)2C \epsilon\}$. Hence, $\norm{\cdot}_{\BB,\Sym,\epsilon}$ is a quasi-norm with modulus of concavity at most $\kappa_\epsilon$. Given $0<q<p$, there is $\epsilon>0$ such that $\kappa_\epsilon=2^{1/q-1}$. Consequently, $\norm{\cdot}_{\BB,\Sym}$ is a quasi-norm equivalent to a $q$-norm. By Lemma~\ref{lem:B0},
\[
\norm{\Ind_{\varepsilon,A}}_\BB \lesssim \abs{A} \lesssim \norm{\Ind_{\varepsilon,A}}_\Sym, \quad A\in\enbrak{N}^{<\infty}.
\]
Hence, \ref{it:CCE} holds. Since $\emptyset\in\Gr(f)$ for all $f\in\FF^\NN$,
\[
\abs{\Ind_I^*(f)} \le \norm{f}_\BB\le \norm{f}_{\BB\cap\Sym}, \quad f\in\FF^{\NN}.
\]
and \ref{it:BSDual} holds.

Note that \ref{it:CCE} implies that $(\ee_n)_{n=1}^\infty$ is a semi-normalized sequence in $\BB\cap\Sym$. Since $\Sym\subset c_0$ continuously, $\BB\cap\Sym\subset c_0$ continuously. So, to prove that $\BB\cap\Sym$ is a sequence space with the required local convexity, it suffices to see its completeness. To that end, we pick a Cauchy sequence $(f_j)_{j=1}^\infty$. Since $(f_j)_{j=1}^\infty$ is a Cauchy sequence relative to the topology of $\Sym$, there is $f\in\Sym$ such that $\lim_j \norm{f_j-f}_\Sym=0$. In particular, $\lim_j f_j=f$ pointwise. Fix $j\in\NN$. By Lemma~\ref{lem:B2} and Lemma~\ref{lem:Approx},
\begin{align*}
\norm{f_j-f}_{\BB,\Sym}
&\le \max\enbrace{\norm{f_j-f}_\Sym, 6 C \liminf_{k\in\NN} \norm{f_j-f_k}_{\Sym} + \norm{f_j-f_k}_{\BB} }\\
&\le \max\{\kappa, 6 C\} \liminf_{k\in\NN} \norm{f_j-f_k}_{\Sym} + \norm{f_j-f_k}_{\BB}\\
&\le \max\{1+\kappa,1+6 C\} \liminf_{k\in\NN} \norm{f_j-f_k}_{\BB,\Sym}.
\end{align*}
Hence, $\lim_j \norm{f_j-f}_{\BB,\Sym}=0$. Since this implies that $f\in\BB$, we are done.

To prove quasi-grediness, pick $f\in\FF^{\NN}$, $A\in\Gr(f)$, $B\in\Gr(f-S_A(f))$, and $I\in\It$. By Lemma~\ref{lem:GSProj},
\[
\abs{\Ind^{*}_{I\setminus B}\enpar{f-S_A(f)}}=\abs{\Ind^{*}_{I\setminus (A\cup B)}(f)}\le \norm{f}_\BB.
\]
Hence, $\norm{f-S_A(f)}_\BB \le \norm{f}_\BB$.

Lemma~\ref{lem:Schauder} implies \ref{it:PSP}. To prove that $\BB\cap\Sym$ is nonseparable, we pick $h=(a_n)_{n=1}^\infty\in\Sym\setminus \ell_1$. Since $h\in\cc_0$, an application of Lemma~\ref{lem:LeibA} yields $A\in\enbrak{\NN}$ so that $g:=\abs{S_A(h)}$ is Leibnizian and has an admissible sequence $\Jt=(J_k)_{k=1}^\infty$ with $\omega(g,\Jt)>0$. Let $(\rho_j)_{j\in J}$ be as in Lemma~\ref{lem:B77}. Set $g_j=M_{\rho_j}(g)$ for all $j\in J$. By unconditionality, $g_j\in\Sym$and $\norm{g_j}_\Sym \le \norm{h}_\Sym$ for all $j\in J$. Since $\norm{g_j}_\BB\le\alpha(g, \Jt)$ and
\[
\norm{g_j-g_k}_{\BB,\Sym} \ge \norm{g_j-g_k}_\BB\ge \omega(g,\Jt)
\]
for all $j\in J$ and $k\in J\setminus\{j\}$, we are done.
\end{proof}

It would be surprising if $\BB\cap\Sym$ improved its convexity relative to $\Sym$. We provide a partial result that points in this direction. A quasi-Banach lattice $\LL$ is said to be \emph{L-convex} if there is $0<\varepsilon<1$ such that
\begin{equation*}
\varepsilon \norm{ f } \le \max_{j\in J} \norm{ f_j}
\end{equation*}
whenever $f\in\LL$ and the finite family $(f_j)_{j\in J}$ in $\LL$ satisfy
\begin{itemize}
\item $(1-\varepsilon) \abs{J} f \le \sum_{j\in J} f_j$, and
\item $0\le f_j\le f$ for every $j\in J$.
\end{itemize}
Most lattices naturally arising analysis are $L$-convex.

\begin{proposition}\label{prop:BNLCG}
Let $0<p< 1$ and $\Sym$ be an $L$-convex symmetric unconditional minimal sequence space with $m \lesssim \usdf[\Sym](m)$ for $m\in\NN$. If $\BB\cap\Sym$ is locally $p$-convex, then $\Sym$ is locally $p$-convex.
\end{proposition}

\begin{proof}
Let $K\in[1,\infty)$ be such that $\norm{f}_\infty\le K \norm{f}_\Sym$ for all $f\in\Sym$. Let $C\in(0,\infty)$ be as in \ref{it:LpC:b} relative to $\norm{\cdot}_{\BB,\Sym}$. Pick a pairwise disjointly supported finite family $(f_j)_{j\in J}$ consisting of nonnegative vectors in $\Sym$. There is a family
\[
(f_{j,n})_{j\in J,n\in\NN}
\]
consisting of nonnegative vectors in $\Sym$ such that $\supp(f_{j,n})\subset \supp{f_j}$ and $f_{j,n}$ is one-to-one on its support for all $(j,n)\in J\times\NN$, and $\lim_n f_{j,n}=f_j$ for all $j\in J$. Since $f_{j,n}\in\cc_0$ for all $j\in J$, for each $n\in\NN$ there is a permutation $\pi_n$ of $\NN$ such that $g_{j,n}:=P_{\pi_n}(f_{j,n})$ is decreasing on its support for all $j\in J$. By Lemma~\ref{lem:LeibB}, for each $j\in J$ there exists a sequence $\Jt_{j,n}=(J_{j,n,k})_{k=1}^\infty$ admissible for $g_{j,n}$ such that $\alpha(g_{j,n},\Jt_{j,n})=\norm{g_{j,n}}_\infty$ and $J_{j,n,k}$ is a singleton contained in $\supp(g_{j,n})$ for all $k\in S(g_{j,n},\Jt_{j,n})$. By Lemma~\ref{lem:B77}, there is $\tau_n\in\{-1,1\}^\NN$ such that $\norm{h_{j,n}}_\BB \le \norm{g_{j,n}}_\infty$ for all $j\in J$, where $h_{j,n}=M_{\tau_n}(g_{j,n})$. We infer that $\norm{h_{j,n}}_{\BB,\Sym}\le K \norm{f_{j,n}}_\Sym$ for all $j\in J$ and $n\in\NN$. Therefore,
\[
\norm{\sum_{j\in J} f_{j,n}}_{\Sym} =\norm{\sum_{j\in J} h_{j,n}}_{\Sym} \le \norm{\sum_{j\in J} h_{j,n}}_{\BB,\Sym}\\
\le C K \enpar{\sum_{j\in J} \norm{f_{j,n}}_{\Sym}^p}^{1/p}.
\]
Let $\kappa$ be the modulus of concavity of $\Sym$. By Lemma~\ref{lem:Approx},
\[
\norm{\sum_{j\in J} f_j}_{\Sym} \le \kappa^2 C K \enpar{\sum_{j\in J} \norm{f_j}_{\Sym}^p}^{1/p}
\]
Applying \cite{Kalton1984b}*{Theorem 2.3} puts an end the proof.
\end{proof}

While the examples of discontinuous quasi-norms that can be found in the literature have been tailored for this purpose, $\norm{\cdot}_{\BB,\Sym}$ is a discontinuous quasi-norm that appears naturally when studying the TGA. We make a detour in our route to substantiate this assertion.

\begin{proposition}\label{prop:Discont}
Let $(\Sym,\norm{\cdot}_\Sym)$ be a minimal unconditional sequence space not contained in $\ell_1$. Assume that $m \lesssim \lsdf[\Sym](m)$ for $m\in\NN$. Then $\norm{\cdot}_{\BB,\Sym}$ restricted to $\cc_{00}$ is a discontinuous map relative to the topology of $\BB\cap\Sym$.
\end{proposition}

\begin{proof}
Let $\kappa\in [1,\infty)$ be such that $\norm{\sum_{j=1}^3 f_j}\le \kappa \sum_{j=1}^3 \norm{f_j}_\Sym$ for all $(f_j)_{j=1}^3$ in $\Sym$. Let $C\in(0,\infty)$ be such that $\norm{f}_\infty\le C\norm{f}_\Sym$ for all $f\in\Sym$. Pick $h\in \cc_{00}$ such that $\norm{h}_\Sym\le 1/C$ and
\[
\norm{h}_1\ge 2+\max\enbrace{1,\kappa \enpar{\norm{\ee_1}_\Sym+\norm{\ee_2}_\Sym+\norm{h}_\Sym}}
\]
We have $\norm{h}_\infty\le 1$. Set $g=S_{\uii{3}}\enpar{\abs{h}}$ and define
\[
G\colon\RR\to \cc_{00}, \quad t\mapsto \ee_1-t \ee_2+ g.
\]
For any $0\le t \le 1$ we have
\begin{multline*}
a:=\norm{g}_1\ge \norm{h}_1 -2\norm{h}_\infty\\
\ge\max\enbrace{1,\kappa \enpar{\norm{\ee_1}_\Sym+\norm{\ee_2}_\Sym+\norm{g}_\Sym}}
\ge \max\enbrace{1, \norm{G(t)}_\Sym}.
\end{multline*}

By Lemma~\ref{lem:BNorm3},
\[
\norm{G(t)}_{\BB,\Sym}
=\begin{cases}
1+a & \mbox{ if } t= 1, \\
1+a -t& \mbox{ if } 0\le t<1.
\end{cases}
\]
Therefore $\norm{G}_{\BB,\Sym}$ is discontinuous at $1$. However, $G\colon\RR\to \BB\cap\Sym$ is continuous.
\end{proof}

Proposition~\ref{prop:BSSS} shows that $\BB\cap\Sym$ may not be minimal. The following result identifies its separable part.

\begin{proposition}
Let $\Sym$ be a minimal sequence space with $m \lesssim \lsdf[\Sym](m)$ for $m\in\NN$. Assume that the unit vector system $\EB$ is a quasi-greedy basis of $\Sym$. Then $\BB_0\cap \Sym$ is the closed subspace of $\BB\cap\Sym$ spanned by the unit vectors.
\end{proposition}

\begin{proof}
The unit vector system of $\BB\cap\Sym$ is quasi-greedy by Proposition~\ref{prop:BSSS}. Hence, by Theorem~\ref{thm:AABWCha}, the separable part of $\BB\cap\Sym$ consists of all $f\in\BB\cap\Sym$ such that
\[
\lim_{A\in\Gr(f)} \norm{f-S_A(f)}_{\BB,\Sym}=0.
\]
Again by Theorem~\ref{thm:AABWCha}, this condition is equivalent to
\[
\lim_{A\in\Gr(f)} \norm{f-S_A(f)}_{\BB}=0.
\]
Since, by Lemma~\ref{lem:GSProj},
\[
\norm{f-S_A(f)}_{\BB}=\sup_{\substack{B\in\Gr(f) \\ B\supset A}} \beta(f,B),
\]
we are done.
\end{proof}

\begin{proposition}\label{prop:BtoBS}
Let $\Sym$ be a minimal sequence space with $m \lesssim \lsdf[\Sym](m)$ for $m\in\NN$. Assume that the unit vector system $\EB$ is a quasi-greedy basis of $\Sym$.
\begin{enumerate}[label=(\roman*),leftmargin=*,widest=iii]
\item\label{it:BtoBS:a} If $\EB$ is a Schauder basis of $\Sym$, so is of $\BB_0\cap \Sym$.
\item\label{it:BtoBS:b} If $\EB$ is superdemocratic relative to $\Sym$, so is relative to $\BB_0\cap \Sym$.
\end{enumerate}
\end{proposition}

\begin{proof}
It is straightforward from Proposition~\ref{prop:BSSS}\ref{it:CCE} and \ref{it:PSP}.
\end{proof}

\begin{corollary}\label{cor:BSAGS}
Let $\Sym$ be a minimal sequence space for which the canonical basis $\EB=(\ee_n)_{n=1}^\infty $
is almost greedy. If $m\lesssim \usdf[\Sym](m)$ for $m\in\NN$, then $\EB$ is an almost greedy basis of $\BB_0\cap \Sym$.
\end{corollary}

\begin{proof}
Just combine Theorem~\ref{thm:CAG} with Proposition~\ref{prop:BtoBS}\ref{it:BtoBS:b}.
\end{proof}

Given a minimal sequence space $\Sym$ such $m \lesssim \lsdf[\Sym](m)$ for $m\in\NN$ and the canonical basis of $\Sym$ is quasi-greedy we denote by
\[
J_{\BB,\Sym}\colon \BB_0\cap \Sym \to \enpar{\BB_0\cap \Sym}^{\wedge}
\]
the Banach envelope map. We point out that $\enpar{\BB_0\cap \Sym}^{\wedge}$ could not be a sequence space. By Proposition~\ref{prop:BtoBS}, a condition that ensures $\enpar{\BB_0\cap \Sym}^{\wedge}$ is a minimal sequence space is that the unit vectors are a Schauder basis of $\Sym$.

To put the following result into context, we recall that $\ell_1$ is the unique symmetric unconditional locally convex sequence space whose fundamental function grows as $(m)_{m=1}^\infty$ (see \cite{LinTza1977}*{Proposition 3.a.6}).

\begin{proposition}\label{prop:EnvDemNUCC}
Let $\Sym$ be a symmetric unconditional minimal sequence space with $\usdf[\Sym](m) \approx m$ for $m\in\NN$. Assume that $\Sym$ is not contained in $\ell_1$. Then the canonical basis of $\enpar{\BB_0\cap \Sym}^{\wedge}$ is democratic and fails to be unconditional for constant coefficients. Moreover,
\[
\norm{\Ind_A}_{\enpar{\BB_0\cap \Sym}^{\wedge}}\approx \abs{A}, \quad A\in\enbrak{\NN}^{<\infty}.
\]
\end{proposition}

\begin{proof}
Set $\tau=((-1)^{n-1})_{n=1}^\infty$. Let $\kappa$ be the modulus of concavity of $\BB\cap\Sym$. Fix $A\in\enbrak{\NN}^{<\infty}$. By Proposition~\ref{prop:BSSS},
\[
\norm{\Ind_A}_{\enpar{\BB_0\cap \Sym}^{\wedge}}
\ge \sup_{I\in\It} \abs{\Ind_I^*(\Ind_A)}= \sup_{I\in\It} \abs{I\cap A}=\abs{A}.
\]
In turn, also by Proposition~\ref{prop:BSSS},
\[
\norm{\Ind_A}_{(\BB_0\cap\Sym)^{\wedge}} \le \norm{\Ind_A}_{\BB,\Sym}\approx \norm{\Ind_A}_{\Sym} \approx \abs{A},
\quad A\in\enbrak{A}^{<\infty}.
\]

Pick a decreasing positive sequence $f\in\Sym\setminus\ell_1$. For $m\in\NN$ put
\[
f_m=S_{\bii{1}{m}}(f).
\]
Let $\pi_m$ be the permutation of $\NN$ given by $\pi_m(n)=n+1$ if $1\le n\le m-1$, $\pi_m(m)=1$, and $\pi(n)=n$ if $n>m$. Note that $\pi_m$ is cyclic on $\bii{1}{m}$. Given $m\in\NN$ and $k\in \bii{0}{m-1}$, the vector
\[
f_{k,m}:=P_{\pi_{m}^k}(f_m)
\]
is decreasing on $\bii{1}{m-k}$ and on $\bii{m-k+1}{m}$. Consequently, by Lemma~\ref{lem:LeibB} and Lemma~\ref{lem:B77},
\[
\norm{M_\tau(f_{k,m})}_{\BB,\Sym}\le D:=2\kappa \max\enbrace{\norm{f}_\Sym,\norm{f}_\infty}.
\]
Therefore, if
\[
g_m=\frac{1}{m}\sum_{k=0}^{m-1} M_\tau(f_{k,m}),
\]
then $\norm{g_m}_{\enpar{\BB_0\cap \Sym}^{\wedge}}\le D$. Set $s_m=\Ind^*_{\bii{1}{m}}(f)$ for all $m\in\NN$. Since
\[
g_m = \frac{s_m}{m} \Ind_{\tau, \bii{1}{m}},
\]
$\norm{\Ind_{\tau, \bii{1}{m}}}_{\enpar{\BB_0\cap \Sym}^{\wedge}}\le Dm/s_m$ for all $m\in\NN$. In turn, if $B_m= 2\ZZ\cap[1,m]$,
\[
\norm{\Ind_{\tau,B_m}}_{(\BB_0\cap\Sym)^{\wedge}}\approx \abs{B_m} \approx m, \quad m\in\NN.
\]
Noticing that $\lim_m s_m=\infty$ puts an end to the proof.
\end{proof}

We conclude this section by gathering the results we have obtained.

\begin{theorem}\label{thm:Bgathers}
There is a quasi-Banach space $\XX$ with an almost greedy Schauder basis $\XB$ whose Banach envelope is democratic and fails to be unconditional for constant coefficients. Moreover,
\[
\norm{\Ind_A[\XB]}\approx \norm{\Ind_A[\widehat{\XB}]}\approx \abs{A}, \quad A\in\enbrak{\NN}^{<\infty},
\]
and $\XX$ is locally $p$-convex for all $0<p<1$.
\end{theorem}

\begin{proof}
Bearing in mind Proposition~\ref{prop:BtoBS}\ref{it:BtoBS:a}, Corollary~\ref{cor:BSAGS} and Proposition~\ref{prop:EnvDemNUCC}, it suffices to show the existence of a symmetric unconditional sequence space $\Sym$ with $\ell_1 \subsetneq \Sym$ whose fundamental function is equivalent to $(m)_{m=1}^\infty$ and which is locally $p$-convex for all $0<p<1$. The Lorentz sequence space $\ell_{1,q}$, $1<q\le\infty$, does the job.
\end{proof}

Note that Theorem~\ref{thm:Bgathers} solves Question~\ref{qt:A} in the negative.
%--------------------------------------------------------------------------
\section{An almost greedy basis that fails to be a Schauder basis under any ordering}\noindent
%--------------------------------------------------------------------------
The construction we carry out in this section is a relative of our construction in Section~\ref{sect:AGDNPE}. Roughly speaking, the difference consists on considering sums $\sum_{n\in\NN\setminus A} a_n$, where $A$ is a greedy set of $(a_n)_{n=1}^\infty$, instead of sums $\sum_{n\in I\setminus A} a_n$, where $I\in\It$. We start by defining a suitable space where these sums are well-defined.

The space we will use instead of $\BB$ is a nonlinear relative of the space of all convergent series,
\[
\st=\enbrace{ (a_n)_{n=1}^\infty \in\FF^\NN \colon \exists \lim_m \sum_{n=1}^m a_n\in \FF}.
\]
Instead of imposing convergence to the partial sums of the series whose general term is $a_n$, $n\in\NN$, we will asume that the series we obtain when arranging $(a_n)_{n=1}^\infty$ acording to the size of $\abs{a_n}$ converges. To be precise, we put
\[
\AB=\enbrace{f\in\cc_0 \colon \enpar{\Ind_A^*(f)}_{A\in\Gr(f)} \mbox{ converges}},
\]
and define
\[
\sigma_g\colon\AB\to \FF, \quad
f\mapsto \lim_{A\in\Gr(f)} \Ind_A^*(f).
\]
In a sense, if $f=(a_n)_{n=1}^\infty\in\AB$, then $\sigma_g(f)$ is the sum of $a_n$ as $n$ runs over $\NN$ \emph{greedily}. We also set
\[
\sigma_g(f,A)= \sigma_g(f)-\Ind_A^*(f), \quad f\in\AB, \, A\in\enbrak{\NN}^{<\infty}.
\]
While $\st$ equipped with the gauge
\[
f=(a_n)_{n=1}^\infty \mapsto \norm{f}_{\st}=\max_{m\in\NN} \abs{\sum_{n=m}^\infty a_n}
\]
is a Banach space isometric to $\cc_0$, we will endow $\AB$ with the gauge
\[
\norm{\cdot}_{\AB} \colon \AB\to [0,\infty), \quad
f\mapsto \norm{f}_{\AB}=\max\enbrace{ \abs{\sigma_g(f,A)} \colon A\in\Gr(f)}.
\]

The following result, which the reader interested exclusively in the answer to Question~\ref{qt:B} may skip, establishes the connections between the spaces $\BB_0$, $\AB$, and $\st$. To state it, for convenience we denote by $\sigma_s$ the summing functional on $\st$, that is,
\[
\sigma_s\colon\st\to \FF, \quad (a_n)_{n=1}^\infty \mapsto \sum_{n=1}^\infty a_n.
\]
We will use the fact that $\st$ is a \emph{spreading} sequence space, that is,
\[
\norm{P_\varphi(f)}_{\st}=\norm{f}_{\st}
\]
for all $f\in\st$ and all increasing maps $\varphi\colon\NN\to\NN$.

\begin{theorem}\label{thm:B0Ac}
$\BB_0\subset \AB \cap \st$. Moreover, if $f \in \BB_0$, then
\[
\max\enbrace{ \norm{f}_{\st},\norm{f}_{\AB}}\le \norm{f}_{\BB},
\]
and $\sigma_g(f)=\sigma_s(f)$. Besides,
\begin{enumerate}[label=(\roman*),leftmargin=*,widest=iii]
\item\label{it:NotEmbedA} $\AB\cap\st\not\subset \BB$,
\item\label{it:NotEmbedB} $\AB\cap\st\cap\BB\not\subset \BB_0$,
\item\label{it:NoSumA} there exist $f_1\in\BB_0$ and $g_1\in\ell_1$ such that $f_1+g_1\notin \BB\cup \AB$, and
\item\label{it:NoSumB} there exist $f_2\in\BB_0$ and $g_2\in\ell_1$ such that $f_2+g_2\in \BB\setminus \AB$.
\end{enumerate}
\end{theorem}

\begin{proof}
Given $A$, $B\in \Gr(f)$, there is $I\in\It$ such that $A\cup B\subset I$. We have
\begin{multline*}
\abs{\Ind_B^*(f)-\Ind_A^*(f)}
=\abs{\Ind_{I\setminus A}^*(f)-\Ind_{I\setminus B}^*(f)}\\
\le\abs{\Ind_{I\setminus A}^*(f)} + \abs{\Ind_{I\setminus B}^*(f)}
\le \beta(f,A)+ \beta(f,B).
\end{multline*}
Therefore, given $D\in\Gr(f)$,
\[
\sup_{\substack{A,B\in\Gr(f) \\ A\cap B\supset D}} \abs{\Ind_B^*(f)-\Ind_A^*(f)}
\le 2 \sup_{\substack{A\in\Gr(f) \\ A\supset D}} \beta(f,A).
\]
Hence, $(\Ind_A^*(f))_{A\in\Gr(f)}$ is a Cauchy net, whence $f\in\AB$. Now, given $A\in\Gr(f)$,
\begin{multline*}
\abs{\sigma_g(f)-\Ind_A^*(f)}
=\lim_{B\in \Gr(f)} \abs{\Ind_B^*(f)-\Ind_A^*(f)} \\
\le \liminf_{B\in\Gr(f)} \beta(f,A)+ \beta(f,B)
=\beta(f,A)
\le \norm{f}_\BB.
\end{multline*}
Therefore, $\norm{f}_\AB\le \norm{f}_\BB$. By Lemma~\ref{lem:GUnion}, for each $k\in\NN$ the set
\[
\enbrace{A\in\Gr(f) \colon A\subset [1,k]}.
\]
has a maximum element, which we denote by $A_k$. If $A\in\Gr(f)\setminus\{\emptyset\}$ and $k\in\uii{\max(A)}$, then $A\subset A_k$. Hence, $(A_k)_{k=1}^\infty$ is a monotone final sequence. Therefore, $\lim_k \beta(f,A_k)=0$.

Set $f=(a_n)_{n=1}^\infty$. If $k\in\NN$ and $m\in\uii{k+1}$, then
\begin{multline*}
\abs{\sum_{n=k+1}^m a_n}
=\abs{\Ind_{\bii{1}{m}\setminus A_k}^*(f)-\Ind_{\bii{1}{k}\setminus A_k}^*(f)}\\
\le \abs{\Ind_{\bii{1}{m}\setminus A_k}^*(f)}+\abs{\Ind_{\bii{1}{k}\setminus A_k}^*(f)}
\le 2 \beta(f,A_k).
\end{multline*}
Consequently, $\sum_{n=1}^\infty a_n$ is a Cauchy series, whence $f\in\st$. Besides, given $k\in\NN$,
\[
\abs{\sum_{n=k}^\infty a_n}=\lim_{m\in\NN} \abs{\Ind^*_{\bii{k}{m}}(f)}\le\norm{f}_\BB.
\]
Therefore, $\norm{f}_{\st}\le\norm{f}_\BB$. Given $A\in\Gr(f)\setminus\{\emptyset\}$,
\begin{equation}\label{eq:GreedyStantard}
\abs{ \sum_{k=1}^{\max(A)} a_n - \Ind^*_A(f)} =\abs {\Ind^*_{\bii{1}{\max(A)} \setminus A}(f)} \le \beta(f,A).
\end{equation}
Since $\abs{A} \le \max(A)$ for all $A\in\enbrak{\NN}\setminus\{\emptyset\}$, $\lim_{A\in\Gr(f)} \max(A)=\infty$. Hence, letting $A$ tend to infinity in \eqref{eq:GreedyStantard} we obtain
\[
\abs{\sigma_s(f)-\sigma_g(f)}=0.
\]

We go on by proving \ref{it:NoSumA} and \ref{it:NoSumB}. Pick a decreasing sequence $(b_k)_{k=0}^\infty$ in $(0,\infty)$, a sequence $(m_k)_{k=1}^\infty$ in $\NN$, and a sequence $(\epsilon_k)_{k=1}^\infty$ in $(0,\infty)$ satisfying the following properties:
\begin{itemize}
\item $\lim_k b_k=0$;
\item when aiming to prove \ref{it:NoSumA} we will asume
\begin{equation}\label{eq:NoSumA}
\sup_{k} m_k b_k=\infty,
\end{equation}
while, when aiming to prove \ref{it:NoSumB},
\begin{equation}\label{eq:NoSumB}
b_{k-1} \approx b_k\approx 1/m_k, \quad k\in\NN;
\end{equation}
\item $\varepsilon_k < b_{k-1}-b_k$ for all $k\in\NN$; and
\item $S:=\sum_{k=1}^\infty m_k \epsilon_k<\infty$.
\end{itemize}

For instance, putting $b_k=2^{-k}$, $m_k=3^k$ (or $m_k=2^k$) and $\epsilon_k=4^{-k}$ serves our purposes. Set
\[
n_k=\sum_{j=1}^k m_j, \quad k\in\uii{0}.
\]
Pick a decreasing sequence $(c_n)_{n=0}^\infty$ such that
\begin{itemize}
\item $c_{2n_k}=b_k$ for all $k\in\uii{0}$,
\item $b_{k-1}-c_{1+2n_{k-1}}>\epsilon_k$ for all $k\in\NN$, and
\item $c_{2+2n_{k-1}}- c_{-1+2n_k}<\epsilon_k$ for all $k\in\NN$.
\end{itemize}
Define $f_0=((-1)^{n-1}c_n)_{n=1}^\infty$ for all $n\in\NN$. By Lemma~\ref{lem:B77}, $f_0\in\BB_0$. Now define $g_0=(d_n)_{n=1}^\infty$ by
\[
d_n=\begin{cases}
\epsilon_k & \mbox{ if $n$ is odd and $k\in\NN$ is such that $n\in \bii{2 n_{k-1}} {2 n_k}$}\\ 0 & \mbox{ if $n$ is even}. \end{cases}
\]
The sequence $g_0$ is nonnegative, and $\sum_{n=1}^\infty d_n=S$.

Set $h_0=f_0+g_0=(r_n)_{n=1}^\infty$. Also set
\[
J_k= \bii{1+2 n_{k-1}} {2 n_k}, \quad J_k^{+}=J_k\cap(1+2\ZZ), \quad J_k^{-}=J_k\cap 2\ZZ,
\]
$I_k =\cup_{j=1}^{k-1} J_j$, and $G_k=I_k \cup J_k^+$ for all $k\in\NN$. We have
\begin{align*}
\max_{n\in J_k^+} \abs{r_n}&=c_{1+2 n_{k-1}} +\epsilon_k< b_{k-1},\\
\min_{n\in J_k^+} \abs{r_n}&=c_{-1+2n_k}+ \epsilon_k> c_{2+2 n_{k-1}},\\
\max_{n\in J_k^-} \abs{r_n}&=c_{2+2 n_{k-1}},\\
\min_{n\in J_k^-} \abs{r_n}&=c_{2 n_{k}} =b_k.
\end{align*}
We infer that
\begin{equation}\label{eq:Bd101}
\abs{\Ind^*_{D}(h_0)} \le 2 m_k b_{k-1}, \quad k\in\NN, \, D\in \enbrak{J_k}.
\end{equation}
We also infer that $A\in\Gr(h_0)$ if and only if there is $k\in\NN$ such that either
$A = I_k \cup B$ for some $B\in\Gr(S_{J_k^+}(h_0))$ or $A = G_k \cup B$ for some $B\in\Gr (S_{J_k^-}(h_0))$. In particular, $\{I_k, G_k\} \subset \Gr(h_0)$ for all $k\in\NN$, and $(G_k)_{k=1}^\infty$ is a monotone final sequence. Since
\[
\abs{\Ind_{I_{k+1}}^*(h_0) - \Ind^*_{G_k}(h_0)} = \abs{\Ind_{J_k\setminus G_k}^*(h_0)}=\abs{\Ind_{J_k^-}^*(h)}\ge m_k b_k
\]
for all $k\in\NN$, $h_0\notin\AB$. Besides, $h_0\notin\BB$ provided that \eqref{eq:NoSumA} holds. Assume, on the contrary, that \eqref{eq:NoSumB} holds. For each $I\in\It$ and $A\in\Gr(h_0)$, there are $j$, $m\in\NN$, $D\in\enbrak{ J_j}$ and $E\in\enbrak{J_{m+1}}$ such that $ j\le m$ and
\[
I\setminus A=D \cup \enpar{\cup_{k=j+1}^m J_k} \cup E.
\]
Since $f_0\in \BB_0\subset \st$ and $g_0\in \ell_1\subset \st$, $h_0\in \st$. Hence, $(\Ind_J^*(h_0))_{J\in\It}$ is a bounded family. Therefore,
\[
\Ind_{I\setminus A}^*(h_0), \quad I\in\It, \, A\in\Gr(h_0),
\]
is a bounded family by \eqref{eq:Bd101}. This means that $h_0\in\BB$. The sequence $h_0$ will also play a key role in our proof of \ref{it:NotEmbedA} and \ref{it:NotEmbedB}. Let $(H_k)_{k=1}^\infty$ be a right-dominant partition of $\NN$ such that $2\abs{J_k}=\abs{H_k}$ for all $k\in\NN$. For each $k\in\NN$, let
$H_k^+$ an $H_k^-$ be the left half and the right half of $H_k$, respectively. Let
\[
\varphi^+\colon\NN\to\cup_{k=1}^\infty H_{k}^+, \quad \varphi^-\colon\NN\to\cup_{k=1}^\infty H_{k}^-
\]
be increasing bijections. Since $\st$ is spreading,
\[
h:=P_{\varphi^+}(h_0)-P_{\varphi^-}(h_0)\in \st.
\]
By \eqref{eq:Bd101},
\begin{equation}\label{eq:Bd103}
\abs{\Ind^*_{D}(h)} \le 4m_k b_{k-1}, \quad k\in\NN, \, D\in \enbrak{H_k}.
\end{equation}

Given $A\in\enbrak{\NN}$, $A\in\Gr(h)$ if and only if $A=\varphi^+(A^+)\cup\varphi^-(A^-)$ for some $A^+$, $A^-\in\Gr(h_0)$ with $\abs{A^+\triangle A^-}\le 1$. If this is the case then $S_A(h)=0$ if $A$ is even, and
\[
\abs{S_A(h)}=D(h_0)\enpar{\frac{\abs{A}+1}{2}}
\]
is $A$ is odd. Hence, $h\in\AB$, and $\sigma_g(h)=0$. Since
\[
H_k^+ \setminus \enpar{\varphi^+(G_k)\cup \varphi^-(G_k)}=\varphi^+(J_k\setminus G_k)
\]
and $\Ind^*_{\varphi^+(J_k\setminus G_k)}(h)=\Ind^*_{J_k\setminus G_k}(h_0)$, $h\notin\BB_0$; and $h\notin\BB$ provided that \eqref{eq:NoSumA} holds. Assume that \eqref{eq:NoSumB} holds and pick $I\in\It$ and $A\in\Gr(h)$. There are $j$, $m\in\NN$, $D\in\enbrak{H_j}$, and $E\in\enbrak{H_{m+1}}$ such that $ j\le m$ and
\[
I\setminus A=D \cup \enpar{\cup_{k=j+1}^m H_k} \cup E.
\]
By \eqref{eq:Bd103}, $h\in\BB$.
\end{proof}

\begin{corollary}
None of the spaces $\BB_0$, $\BB$, $\AB$, $\AB\cap\st$, $\BB\cap\st$ or $\AB\cap\BB\cap\st$ are linear.
\end{corollary}

\begin{proof}
Just combine Theorem~\ref{thm:B0Ac}\ref{it:NoSumA} with Lemma~\ref{lem:B0}.
\end{proof}

From now on, we will focus on building the whised-for space that answers Question~\ref{qt:B}.

\begin{lemma}\label{lem:AConv}
Given $f\in\AB$ and $A\in\Gr(f)$, $f-S_A(f)\in\AB$ and
\[
\sigma_g\enpar{f-S_A(f)}=\sigma_g(f,A).
\]
Moreover, $\lim_{A\in\Gr(f)} \norm{f-S_A(f)}_{\AB}=0$.
\end{lemma}

\begin{proof}
Set $g:=f-S_A(f)\in\AB$. For any $B\in\enbrak{\NN}^{<\infty}$ we have
\[
\sigma_g(f,A)-\Ind_B^*(g)=\sigma_g(f,A\cup B).
\]
By Lemma~\ref{lem:GSProj}, the mapping $B\mapsto A\cup B$ defines a monotone final function from $\Gr(g)$ to $\Gr(f)$. Therefore,
\[
\lim_{B\in\Gr(g)} \sigma_g(f,A\cup B)=0.
\]
This means that $g\in\AB$ and $\sigma_g(g)=\sigma_g(f,A)$. Besides,
\[
\norm{g}_{\AB}\le M(A):=\sup\enbrace{ \abs{\sigma_g(f,D)} \colon D\in\Gr(f), \, D\supset A}.
\]
Since $\lim_{A\in\Gr(f)} M(A)=0$, we are done.
\end{proof}

\begin{lemma}[cf.\@ Lemma~\ref{lem:B1}]
Let $f$, $g\in h_{1,\infty} \cap \AB$. Then, $f+g\in \AB$, $\sigma_g(f+g)=\sigma_g(f)+\sigma_g(g)$, and $\norm{f+g}_\AB\le N(f,g)$, where
\[
N(f,g):= 2 \enpar{ \norm{f+g}_{1,\infty}+\norm{f}_{1,\infty}+\norm{g}_{1,\infty}}+\norm{f}_\AB+\norm{g}_\AB.
\]
\end{lemma}

\begin{proof}
Pick $A\in\Gr(f+g)$. Choose $B\in\Gr(f)$ and $D\in\Gr(D)$ with $m:=\abs{A}=\abs{B}=\abs{D}$. Choose $k\in\bii{0}{m-1}$. By \eqref{eq:twisted} and Lemma~\ref{lem:B99},
\begin{align*}
N(A):=&\abs{\sigma_g(f)+\sigma_g(g)-\Ind^*_A(f+g)} \le \abs{\sigma_g(f,B)}+\abs{\sigma_g(g,D)} \\
&+ \frac{2 m}{1+m-k} \enpar{\rho_{1,\infty}(f,k)+ \rho_{1,\infty}(g,k) + \rho_{1,\infty}(f+g,k) }.
\end{align*}
By choosing $k=\enfloor{m/2}$ we get $\lim_{A\in\Gr(f+g)} N(A)=0$, while if we choose  $k=0$ we obtain $N(A)\le N(f,g)$.
\end{proof}

The proofs of the following two results go over the lines of the corresponding results for $\norm{\cdot}_{\BB}$, thus we will omit their proofs.

\begin{lemma}[cf.\@ Lemma~\ref{lem:BNorm3}]
Let $(a_n)_{n=1}^\infty\in \cc_{00}$. Assume that $a_n\in[0,1]$ for all $n\in\NN$, $a_1=a_2=0$, and $a:=\sum_{n=1}^\infty a_n \ge 1$. Then,
\[
\norm{\ee_1- t \ee_2+f}_\AB
=\begin{cases}
1+a & \mbox{ if } t= 1, \\
1+a -t& \mbox{ if } 0\le t<1.
\end{cases}
\]
\end{lemma}

\begin{lemma}[cf.\@ Lemma~\ref{lem:B2}]\label{lem:A2}
Given $f\in\FF^\NN$ and $(f_j)_{j=1}^\infty$ in $\AB$ such that $\lim_j f_j=f$ pointwise,
\[
\sup\enbrace{\abs{\Ind^*_A(f)} \colon A\in\Gr(f)} \le \liminf_{j\in\NN} 6\norm{f_j}_{1,\infty} + 2 \norm{f_j}_\AB.
\]
Moreover, if $f\in\AB$, then
\[
\norm{f}_{\AB} \le \liminf_{j\in\NN} 6 \norm{f_j}_{1,\infty} + \norm{f_j}_\AB.
\]
\end{lemma}

Given a sequence space $(\Sym,\norm{\cdot}_\Sym)$ we set
\[
\norm{f}_{\AB,\Sym}=\max\enbrace{ \norm{f}_\Sym, \norm{f}_\AB}, \quad f\in\AB\cap\Sym.
\]

In the following result, we collect some relevant properties of $\AB\cap\Sym$. Before doing so, we make some comments that may help the reader understand the similarities and dissimilarities between $\AB\cap\Sym$ and $\BB\cap\Sym$. On  one hand, Lemma~\ref{lem:Schauder} no longer holds for $\norm{\cdot}_{\AB}$. On the other hand $\norm{\cdot}_{\AB}$, unlike $\norm{\cdot}_{\BB}$, is symmetric. These features of $\norm{\cdot}_{\AB}$ will cause the minimal systems built from it to be symmetric, but not Schauder bases. Conversely, the minimal systems built from $\norm{\cdot}_{\BB}$ are non-symmetric Schauder bases.

\begin{proposition}[cf.\@ Propositions~\ref{prop:BtoBS}, \ref{prop:EnvDemNUCC} and Proposition~\ref{prop:Discont}]\label{prop:ASSS}
Let $(\Sym,\norm{\cdot}_\Sym)$ be a minimal sequence space. Suppose the canonical basis of $\Sym$ is quasi-greedy with $m \lesssim \lsdf[\Sym](m)$ for $m\in\NN$. Then $(\AB\cap\Sym,\norm{\cdot}_{\AB,\Sym})$ is a minimal sequence space whose canonical basis is quasi-greedy. Furthermore,
\begin{enumerate}[label=(\roman*),leftmargin=*,widest=iii]
\item $\norm{\Ind_{\varepsilon, A}}_{\AB,\Sym} \approx \norm{\Ind_{\varepsilon, A}}_\Sym$ for $A\in\enbrak{\NN}^{<\infty}$ and $\varepsilon\in S_\FF^A$.
\item $\sigma_g$ is a linear operator on $\AB\cap\Sym$, and $\norm{\sigma_g}_{\AB\cap\Sym\to\FF} \le 1$.
\item If $\Sym$ is locally $p$-convex, $0<p\le 1$, then $\AB\cap\Sym$ is locally $q$-convex for all $0<q<p$.
\item If $\Sym$ is symmetric, so is $\AB\cap\Sym$.
\item If the unit vector system $\EB$ is superdemocratic relative to $\Sym$, so is relative to $\AB\cap\Sym$.
\item If $\EB$ is an almost greedy basis of $\Sym$, so is of $\AB\cap\Sym$.
\item If $\EB$ is an unconditional basis of $\Sym$ and $\Sym$ is not contained in $\ell_1$, then $\norm{\cdot}_{\AB,\Sym}$ restricted to $\cc_{00}$ is a discontinuous map relative to the topology of $\AB\cap\Sym$.
\end{enumerate}
\end{proposition}

\begin{proof}
We will focus on the points where the arguments diverge from those of the proof of Proposition~\ref{prop:BtoBS}, and omit other details. To prove completeness, we pick  a Cauchy sequence $(f_j)_{j=1}^\infty$ in $\AB\cap\Sym$. Put
\[
\epsilon_j=\limsup_{k\in\NN} \norm{f_j-f_k}_\AB, \quad j\in\NN.
\]
There are $f\in\Sym$ with $\lim_j f_j =f$ (in the topology of $\Sym$) and $s\in\RR$ such that $\lim_j \sigma_g(f_j)=s$. Let $\gamma_j=\abs{s-\sigma_g(f_j)}$, $g_j=f_j-f$, and $\delta_j= \norm{g_j}_{1,\infty}$ for all $j\in\NN$. Also set $\lambda_m= \rho_{1,\infty}\enpar{f,\enfloor{m/2}}$ for all $m\in\NN$. By Theorem~\ref{thm:LorentzEmbeds},
\[
\lim_m \lambda_m=\lim_{j\in\NN} \epsilon_j=\lim_{j\in\NN} \delta_j= \lim_{j\in\NN} \gamma_j=0.
\]

We will use twice the following consequence of Lemma~\ref{lem:A2}.
\begin{claim}\label{claim:A}
For any $j\in\NN$ and $D\in\Gr(g_j)$, $\abs{\Ind^*_D(g_j)} \le 6 \delta_j +2 \epsilon_j$. Moreover, if $f$ belongs to $\AB$ then $\norm{g_j}_\AB\le 6\delta_j + \epsilon_j$.
\end{claim}

Set $\lambda_m= \rho_{1,\infty}\enpar{f,\enfloor{m/2}}$ for all $m\in\NN$. Combining inequality~\eqref{eq:twisted}, Lemma~\ref{lem:B99} and Lemma~\ref{lem:BS} gives the following.

\begin{claim}\label{claim:B}
If $A\in\Gr(f)$, $B\in\Gr(f_j)$, and $D\in\Gr(g_j)$ are such that $m=\abs{A}=\abs{B}=\abs{D}$, then
\[
\abs{\Ind^*_B(f_j)-\Ind^*_A(f)-\Ind^*_D(g_j)} \le 10 \delta_j +12 \lambda_m.
\]
\end{claim}
Indeed, if $k=\enfloor{m/2}$, the left-hand term of the inequality is bounded above by
\begin{align*}
&\frac{2 m}{1+m} \norm{g_j}_{1,\infty}+ \frac{2 m}{m-k} \rho_{1,\infty}(f_j,1+k)+\frac{2 m}{1+m-k} \rho_{1,\infty}(f,k)\\
&\le \enpar{ \frac{2 m}{1+m}+\frac{4m}{m-k}}\norm{g_j}_{1,\infty}
+\enpar{\frac{4m}{m-k} + \frac{2 m}{1+m-k} } \rho_{1,\infty}(f,k).
\end{align*}

We return to the main line of the proof. For each $A\in\Gr(f)$ and each $j\in\NN$ pick $B_j=B_j(A)\in\Gr(f_j)$ and $D_j=D_j(A)\in\Gr(g_j)$ with $\abs{D_j}=\abs{B_j}=\abs{A}$. By Claim~\ref{claim:A} and Claim~\ref{claim:B}, $\abs{s-\sigma_g(f,A)}$ is bounded above by
\begin{multline*}
\gamma_j+\abs{\sigma_g(f_j)-\Ind_{B_j}^*(f_j)}
+\abs{\Ind^*_{B_j}(f_j)-\Ind^*_A(f)-\Ind^*_{D_j}(g_j)}+\abs{\Ind^*_{D_j}(g_j)}\\
\le \abs{\sigma_g(f_j)-\Ind_{B_j}^*(f_j)} + 2 \epsilon_j+ \gamma_j + 16 \delta_j +12 \lambda_m.
\end{multline*}
Since $\lim_{A\in\Gr(f)} \lambda_{\abs{A}}=0$ and, for any fixed $j\in\NN$,
\[
\lim_{A\in\Gr(f)} \abs{\sigma_g(f_j)-\Ind_{B_j}^*(f_j)}=0,
\]
$\lim_{A\in\Gr(f)} \abs{s-\sigma_g(f,A)}=0$. This proves that $f\in\AB$ and $\sigma_g(f)=s$. By Claim~\ref{claim:A}, $\lim_j \norm{f-f_j}_{\AB,\Sym}=0$.

Finally we prove that $\cc_{00}$ is dense in $\AB\cap\Sym$. Pick $f\in\AB\cap\Sym$. By Theorem~\ref{thm:AABWCha} and Lemma~\ref{lem:AConv},
\[
\lim_{A\in\Gr(f)} \norm{f-S_A(f)}_{\AB,\Sym}=0.\qedhere
\]
\end{proof}

Given sequence spaces $(\XX, \norm{\cdot}_\XX)$ and $(\YY,\norm{\cdot}_\YY)$, then $\XX\cap \YY$ endowed with the gauge
\[
f\mapsto \max\enbrace{ \norm{f}_\XX, \norm{f}_\YY}
\]
is a sequence space.

\begin{corollary}
Let $\Sym$ be a minimal sequence space with $m \lesssim \lsdf[\Sym](m)$ for $m\in\NN$. Assume that the canonical basis of $\Sym$ is quasi-greedy. Then $\BB_0\cap \Sym\subset \enpar{\AB\cap\Sym}\cap \st$ and the inclusion is a contraction.
\end{corollary}

\begin{proof}
It is a ready consequence of Theorem~\ref{thm:B0Ac}.
\end{proof}

The proofs of the following two results also go over the lines of the corresponding results for $\norm{\cdot}_{\BB}$. We are mainly interested in the latter, and include the former for the sake of completeness.

\begin{proposition}[cf.\@ Proposition~\ref{prop:BNLCG}]
Let $0<p< 1$ and $\Sym$ be an L-convex symmetric unconditional minimal sequence space with $m \lesssim \usdf[\Sym](m)$ for $m\in\NN$. If $\AB\cap\Sym$ is locally $p$-convex, then $\Sym$ is locally $p$-convex.
\end{proposition}

\begin{remark}
In a sense, Proposition~\ref{prop:ASSS} evinces that $\AB\cap h_{1,\infty}$ is the optimal linear space for defining sums of series rearranged according to a greedy criterion. As we have seen, this space fails to be locally convex. However, $\AB\cap h_{1,\infty}$ is locally $p$-convex for all $0<p<1$.
\end{remark}

Given a minimal sequence space $\Sym$ whose canonical basis is quasi-greedy and $m \lesssim \lsdf[\Sym](m)$ for $m\in\NN$, we denote by
\[
J_{\AB,\Sym} \colon \AB\cap\Sym \to \enpar{\AB\cap\Sym}^{\wedge}
\]
the Banach envelope map.

\begin{proposition}[cf.\@ Proposition~\ref{prop:EnvDemNUCC}]\label{prop:AEnvDemNUCC}
Let $\Sym$ be a symmetric unconditional minimal sequence space with $\usdf[\Sym](m) \approx m$ for $m\in\NN$. Assume that $\Sym$ is not contained in $\ell_1$. Then the canonical basis of $(\AB\cap\Sym)^{\wedge}$ is democratic and fails to be unconditional for constant coefficients. Moreover,
\[
\norm{\Ind_A}_{(\AB\cap\Sym)^{\wedge}}\approx \abs{A}, \quad A\in\enbrak{\NN}^{<\infty}.
\]
\end{proposition}

Unlike $\enpar{\BB_0\cap\Sym}^{\wedge}$, $\enpar{\AB\cap\Sym}^{\wedge}$ may fail to be a sequence space.

\begin{proposition}\label{prop:NSS}
Let $\Sym$ be a symmetric unconditional minimal sequence space with $m \lesssim \lsdf[\Sym](m)$ for $m\in\NN$. Assume that $\Sym$ is not contained in $\ell_1$. Then $\widehat{\EB}:=\enpar{J_{\AB,\Sym}(\ee_n)}_{n=1}^\infty$ fails to be total.
\end{proposition}

\begin{proof}
Note that $\widehat{\EB}$ is symmetric. If $\widehat{\EB}$ were total, then it would be unconditional by Theorem~\ref{thm:symunc}. In particular, $\widehat{\EB}$ would be unconditional for constant coefficients. This absurdity proves the result.
\end{proof}

The following result concludes our analysis  of $\AB\cap\Sym$ in the same way as Theorem~\ref{thm:Bgathers} culminates our analysis of $\BB\cap\Sym$.

\begin{theorem}\label{thm:Agathers}
There is a quasi-Banach space $\YY$ with a symmetric almost greedy basis $\YB$ whose envelope $\widehat{\YB}$ is democratic, fails to be total, and fails to be unconditional for constant coefficients. Moreover,
\[
\norm{\Ind_A[\YB]}\approx \norm{\Ind_A[\widehat{\YB}]}\approx \abs{A}, \quad A\in\enbrak{\NN}^{<\infty},
\]
and $\YY$ is locally $p$-convex for all $0<p<1$.
\end{theorem}

\begin{proof}
Apply Proposition~\ref{prop:ASSS}, Proposition~\ref{prop:AEnvDemNUCC}, and Proposition~\ref{prop:NSS} with $\Sym=\ell_{1,q}$, $1<q<\infty$, or $\Sym=h_{1,\infty}$.
\end{proof}

We are now ready to answer Question~\ref{qt:B} in the negative.

\begin{theorem}\label{thm:AMain}
There is a quasi-Banach space $\YY$ with an almost greedy basis $\YB$ that fails to be a Schauder basis under any rearrangement. Furthermore,  $\YB$ can be chosen to be symmetric.
\end{theorem}

\begin{proof}
Let $\YY$ and $\YB$ be as in Theorem~\ref{thm:Agathers}. If $\YB$ were a Schauder basis, then $\widehat{\YB}$ would be a Schauder basis. In particular, $\widehat{\YB}$ would be total. This absurdity proves the result.
\end{proof}

We emphasize that our proof of Theorem~\ref{thm:AMain} takes advantage of the idiosyncracies  of nonlocally convex spaces. In fact, we should concede that Question~\ref{qt:B} remains open when we restrict ourselves to Banach spaces.

As a by-product of our study we show that Theorem~\ref{thm:symunc} breaks down in the quasi-Banach setting.

\begin{corollary}
There is a quasi-Banach space with a symmetric Markushevich basis that fails to be unconditional.
\end{corollary}

\begin{proof}
The Markushevich basis provided by Theorem~\ref{thm:AMain} fails to be unconditional.
\end{proof}
%--------------------------------------------------------------------------
\subsection*{Conflict of Interest}\noindent
%--------------------------------------------------------------------------
The authors declare that they have no conflict of interest.
% ------------------------------------------------------------------------

% ------------------------------------------------------------------------

\begin{bibdiv}
\begin{biblist}

\bib{AABW2021}{article}{
author={Albiac, Fernando},
author={Ansorena, Jos\'{e}~L.},
author={Bern\'{a}, Pablo~M.},
author={Wojtaszczyk, Przemys{\l}aw},
title={Greedy approximation for biorthogonal systems in quasi-{B}anach spaces},
date={2021},
journal={Dissertationes Math. (Rozprawy Mat.)},
volume={560},
pages={1\ndash 88},
}

\bib{AABCO2024}{article}{
author={Albiac, Fernando},
author={Ansorena, Jos\'{e}~L.},
author={Blasco, \'{O}scar},
author={Chu, H\`{u}ng~Viet},
author={Oikhberg, Timur},
title={Counterexamples in isometric theory of symmetric and greedy bases},
date={2024},
ISSN={0021-9045,1096-0430},
journal={J. Approx. Theory},
volume={297},
pages={Paper No. 105970, 20},
url={https://doi.org/10.1016/j.jat.2023.105970},
review={\MR{4650744}},
}

\bib{AACD2018}{article}{
author={Albiac, Fernando},
author={Ansorena, Jos\'{e}~L.},
author={C\'{u}th, Marek},
author={Doucha, Michal},
title={Lipschitz free {$p$}-spaces for {$0 < p < 1$}},
date={2020},
ISSN={0021-2172},
journal={Israel J. Math.},
volume={240},
number={1},
pages={65\ndash 98},
url={https://doi-org/10.1007/s11856-020-2061-5},
review={\MR{4193127}},
}

\bib{AAT2025}{article}{
author={Albiac, Fernando},
author={Ansorena, Jose~L.},
author={Temlyakov, Vladimir},
title={Twenty-five years of greedy bases},
date={2025},
ISSN={0021-9045},
journal={Journal of Approximation Theory},
volume={307},
pages={106141},
}

\bib{AlbiacKalton2016}{book}{
author={Albiac, Fernando},
author={Kalton, Nigel~J.},
title={Topics in {B}anach space theory},
edition={Second},
series={Graduate Texts in Mathematics},
publisher={Springer, [Cham]},
date={2016},
volume={233},
ISBN={978-3-319-31555-3; 978-3-319-31557-7},
url={https://doi.org/10.1007/978-3-319-31557-7},
note={With a foreword by Gilles Godefroy},
review={\MR{3526021}},
}

\bib{Aoki1942}{article}{
author={Aoki, Tosio},
title={Locally bounded linear topological spaces},
date={1942},
ISSN={0369-9846},
journal={Proc. Imp. Acad. Tokyo},
volume={18},
pages={588\ndash 594},
url={http://projecteuclid.org/euclid.pja/1195573733},
review={\MR{14182}},
}

\bib{CRS2007}{article}{
author={Carro, Mar\'{\i}a~J.},
author={Raposo, Jos\'{e}~A.},
author={Soria, Javier},
title={Recent developments in the theory of {L}orentz spaces and weighted inequalities},
date={2007},
ISSN={0065-9266},
journal={Mem. Amer. Math. Soc.},
volume={187},
number={877},
pages={xii+128},
url={https://doi-org/10.1090/memo/0877},
review={\MR{2308059}},
}

\bib{DKKT2003}{article}{
author={Dilworth, Stephen~J.},
author={Kalton, Nigel~J.},
author={Kutzarova, Denka},
author={Temlyakov, Vladimir~N.},
title={The thresholding greedy algorithm, greedy bases, and duality},
date={2003},
ISSN={0176-4276},
journal={Constr. Approx.},
volume={19},
number={4},
pages={575\ndash 597},
url={https://doi-org/10.1007/s00365-002-0525-y},
review={\MR{1998906}},
}

\bib{Hunt1966}{article}{
author={Hunt, R.~A.},
title={On {$L(p,\,q)$} spaces},
date={1966},
ISSN={0013-8584},
journal={Enseign. Math. (2)},
volume={12},
pages={249\ndash 276},
review={\MR{223874}},
}

\bib{Hyers1939}{article}{
author={Hyers, D.~H.},
title={Locally bounded linear topological spaces},
date={1939},
ISSN={0034-7760},
journal={Rev. Ci. (Lima)},
volume={41},
pages={555\ndash 574},
review={\MR{1915}},
}

\bib{Kalton1984b}{article}{
author={Kalton, Nigel~J.},
title={Convexity conditions for nonlocally convex lattices},
date={1984},
ISSN={0017-0895},
journal={Glasgow Math. J.},
volume={25},
number={2},
pages={141\ndash 152},
url={https://doi-org/10.1017/S0017089500005553},
review={\MR{752808}},
}

\bib{KPR1984}{book}{
author={Kalton, Nigel~J.},
author={Peck, N.~Tenney},
author={Roberts, James~W.},
title={An {$F$}-space sampler},
series={London Mathematical Society Lecture Note Series},
publisher={Cambridge University Press, Cambridge},
date={1984},
volume={89},
ISBN={0-521-27585-7},
url={https://doi.org/10.1017/CBO9780511662447},
review={\MR{808777}},
}

\bib{KoTe1999}{article}{
author={Konyagin, Sergei~V.},
author={Temlyakov, Vladimir~N.},
title={A remark on greedy approximation in {B}anach spaces},
date={1999},
ISSN={1310-6236},
journal={East J. Approx.},
volume={5},
number={3},
pages={365\ndash 379},
review={\MR{1716087}},
}

\bib{LinTza1977}{book}{
author={Lindenstrauss, Joram},
author={Tzafriri, Lior},
title={Classical {B}anach spaces. {I} -- sequence spaces},
series={Ergebnisse der Mathematik und ihrer Grenzgebiete [Results in Mathematics and Related Areas]},
publisher={Springer-Verlag, Berlin-New York},
date={1977},
ISBN={3-540-08072-4},
review={\MR{0500056}},
}

\bib{Rolewicz1957}{article}{
author={Rolewicz, Stefan},
title={On a certain class of linear metric spaces},
date={1957},
journal={Bull. Acad. Polon. Sci. Cl. III.},
volume={5},
pages={471\ndash 473, XL},
review={\MR{0088682}},
}

\bib{Woj2000}{article}{
author={Wojtaszczyk, Przemys{\l}aw},
title={Greedy algorithm for general biorthogonal systems},
date={2000},
ISSN={0021-9045},
journal={J. Approx. Theory},
volume={107},
number={2},
pages={293\ndash 314},
url={https://doi-org/10.1006/jath.2000.3512},
review={\MR{1806955}},
}

\end{biblist}
\end{bibdiv}
\end{document}